\documentclass[preprint,11pt]{elsarticle}

\usepackage[a4paper, total={16cm, 22cm}]{geometry}
\usepackage{enumitem}
\usepackage[english]{babel}
\usepackage{xcolor}
\usepackage[colorinlistoftodos,prependcaption]{todonotes}
\usepackage{subcaption}
\usepackage[linesnumbered,ruled,vlined]{algorithm2e}
\usepackage{graphicx}
\usepackage{tikz}
\usetikzlibrary{patterns}
\usepackage{tikz-3dplot}
\usepackage{pgfplots}
\usepgfplotslibrary{fillbetween}
\pgfplotsset{compat=1.3}
\usepackage[utf8]{inputenc}
\usepackage{amsthm,amsmath,amssymb,amsfonts,eurosym}
\newtheorem{theorem}{Theorem}

\newtheorem{definition}{Definition}

\newtheorem{proposition}{Proposition}

\biboptions{round,authoryear,longnamesfirst}
\newcommand\Tstrut{\rule{0pt}{2.6ex}}       
\newcommand\Bstrut{\rule[-0.9ex]{0pt}{0pt}} 
\newcommand{\TBstrut}{\Tstrut\Bstrut} 

\makeatletter

\def\ps@pprintTitle{%

  \let\@oddhead\@empty

  \let\@evenhead\@empty

  \def\@oddfoot{\reset@font\hfil\thepage\hfil}

  \let\@evenfoot\@oddfoot

}

\makeatother

\begin{document}

\begin{frontmatter}

\title{New $\epsilon-$constraint methods for multi-objective integer linear programming: a Pareto front representation approach}                      
\author[cegist]{Mariana Mesquita-Cunha\corref{cor1}}
\ead{mariana.cunha@tecnico.ulisboa.pt}

\author[cegist]{José Rui Figueira}

\author[cegist]{Ana Paula Barbosa-Póvoa}

\address[cegist]{CEGIST - Centre for Management Studies, Instituto Superior Técnico, Universidade de Lisboa, Portugal}

\cortext[cor1]{Corresponding author}

\begin{abstract}
	Dealing with multi-objective problems by using generation methods has some interesting advantages since it provides the decision-maker with the complete information about the set of non-dominated points (Pareto front) and a clear overview of the problem. However, providing many solutions to the decision-maker might also be overwhelming. As an alternative approach, presenting a representative set of solutions of the Pareto front may be advantageous. Choosing such a representative set is by itself also a multi-objective problem that must consider the number of solutions to present, the uniformity, and/or the coverage of the representation, to guarantee its quality. This paper proposes three algorithms for the representation problem for multi-objective integer linear programming problems with two or more objective functions, each one of them dealing with each dimension of the problem (cardinality, coverage, and uniformity). Such algorithms are all based on the $\epsilon$-constraint approach. In addition, the paper also presents strategies to overcome poor estimations of the Pareto front bounds. The algorithms were tested on the ability to efficiently generate the whole Pareto front or its representation. The uniformity and cardinality algorithms proved to be very efficient both in binary and integer problems, being amongst the best in the literature. Both coverage and uniformity algorithms provide good quality representations on their targeted objective, while the cardinality algorithm appears to be the most flexible, privileging uniformity for lower cardinality representations and coverage on higher cardinality. 
\end{abstract}

\end{frontmatter}



\section{Introduction}

\noindent Most real-world problems, although being multi-objective by their very nature, are frequently modelled as single objective problems mainly due to the the lack of multi-objective tools and/or the computational complexity of solving multi-objective models. However, the manipulation of objectives in order to amalgamate them into a single objective function tends to oversimplify the problem and the conflicting nature among objectives. Indeed, a single objective model may not represent appropriately the decision-maker's (DM) real preferences and objectives \citep{BrankeEtAl2008}. The aforementioned disadvantages together with the rapid rise of computational performance, both in terms of hardware and software, are making multi-objective models, and more specifically Multi-Objective Optimization (MOO), a more favoured approach. Nevertheless, multi-objective integer and mixed integer linear programming (MOILP and MOMILP), that arise in many real-world applications (such as logistics and production problems), have received less attention when compared to binary and continuous problems \citep{AlvesClimaco2007}.

Multi-objective models can be classified, depending on the stage at which the DM is involved to ascertain preferences between solutions, into \textit{a priori} methods, interactive methods, and \textit{a posteriori}/generation methods. In \textit{a priori} methods the DM intervenes at the very beginning, stating her/his preferences and the problem is then solved by aggregating the objective functions. This approach not only carries much of the disadvantages presented above for the single objective models, but also most of the times, in complex functions, the DM struggles to provide parameters, as for example, the weights for each objective. In interactive methods, the DM iteratively states and adjusts preferences based on the results previously obtained. The drawback with this approach is two-folded: (1) it may take a long time for a satisfactory trade-off between objectives to be found; and (2) the DM, by never seeing the full Pareto front, also never clearly understands the relation between objectives and their impact on the solutions. At last, \textit{a posteriori} or generation methods, which involve the DM only after finding the Pareto front, are the most advantageous. This approach clearly allows for defining the trade-off between all objectives, providing the DM a full view of the problem and all the tools for making a better informed decision. However, it may become not only more computationally heavy, but also risks being overwhelming for a DM that would not require so many alternatives.

The \textit{a posteriori} approaches require significant computational effort and time to compute the full efficient/non-dominated set. Hence, apart from usually small and linear problems, the full non-dominated set is not practical to compute in a reasonable amount of time. As a result, there are two classes of generation methods, one that aims to determine the whole set of non-dominated solutions (the so called Pareto front), and the other  which focuses on obtaining a set of solutions representative of the non-dominated set without being overwhelming for the DM \citep{AlvesClimaco2007,KiddEtAl2020}. To address the latter, the representation problem must be considered, which, in itself is a multi-objective problem with three objective functions: cardinality, the number of solutions presented to the DM; coverage, how well the complete Pareto frontier is being represented by the set of solutions; and uniformity, how well those solutions are spread through the Pareto front space \citep{Sayin2000}.

In this work we propose three multi-objective \textit{a posteriori} algorithms to solve MOILP/ MOMILP with more than two objectives. Common to the three algorithms, their structure presents strategies to overcome poor quality approximations for the bound of the Pareto front, and can be used either to generate the complete Pareto front or to compute a representation of it. To this end, the algorithms we put forward tackle each of the three dimensions of the representation problem.


Using a generation method on a MOILP/MOMILP is often more challenging than on a multi-objective linear programming since the former has a non-convex feasible region, which implies that there may exist unsupported non-dominated solutions. \cite{ChalmetEtAl1986} proposed a modified weighted-sum method able to produce, apart from the set of supported non-dominated solutions, the whole set unsupported non-dominated solutions. This was achieved by adding constraints to the problem that bound the values of the objective functions. Although, all non-dominated solutions can be obtained with a full parametrization of the weights added to each objective function, it can lead to an extensive and computationally demanding optimization problems.

Several other authors also proposed approaches based on a sequential reduction of the feasible region \citep{KleinHannan1982,SylvaCrema2004,SylvaCrema2007}. \cite{KleinHannan1982} proposes a reduction of the feasible region by sequentially adding constraints that eliminate solutions dominated by the previously found non-dominated solution. This strategy finds solutions spaced by at least a fixed amount (provided by a parameter), for each objective. However, it may skip interesting solutions from the DM preferences point of view, since only one of the objective is considered as the objective function to be optimized. \cite{SylvaCrema2004} presented a variation of \cite{KleinHannan1982}'s method where all the objectives are included in the objective function by using weighting parameters. \cite{SylvaCrema2007} builds upon the \cite{SylvaCrema2004} method by determining, at each iteration, the weights that maximizes the infinity-norm distance to the set dominated by the previously found solutions. Consequently, uniformly spaced solutions are found. All of these methods share the same main drawback, namely the fact that the problem size increases whenever a non-dominated solution is found, since  the new constraints and, in some cases, new variables have to be added to the model. \cite{KiddEtAl2020} also presented a scalarization algorithm but targeted for generating representations, as opposed to the complete Pareto front, for bi-objective problems. To that end, an insertion method based on Voronoi cuts used to partition the search region was developed. This method was proven to achieve simultaneously good coverage and uniformity for a given cardinality. The algorithm continues to insert solutions until a desired cardinality level is reached.  


One of the most popular methods for solving \textit{a posteriori} MOILP problems is the $\epsilon$-constraint method. For constraining each objective, the method resorts to a virtual grid spaced, for each objective, by making use of an $\epsilon$ parameter. \cite{LaumannsEtAl2006} applied the $\epsilon$-constraint method as a generation method by dividing, after each iteration, the objective space using the values of the previous computed solution and limiting the search space sequentially. Additionally, after each new solution is found, a lexicographic optimization is performed to ensure getting non-dominated solutions. \cite{KirlikSayin2014} developed a method based on the \cite{LaumannsEtAl2006} works, which makes use of a two stage model in order to guarantee non-dominated solutions. It is  a search approach based on the construction of rectangles, as in \cite{LaumannsEtAl2006}, but it also removed rectangles in which no non-dominated solution can be found. This approach makes the search less exhaustive when compared to the \cite{LaumannsEtAl2006}'s algorithm. \cite{Mavrotas2009} adapted the original $\epsilon$-constraint method by introducing slack variables in equality constraints and incorporating them in the objective function weighted by a factor which includes the range of each objective. To choose the $\epsilon$ vector, the range of each objective is divided into evenly spaced intervals generating a uniform grid. Nonetheless, this approach produces several redundant solutions. To overcome such a problem, \cite{MavrotasFlorios2013} extended the  \cite{Mavrotas2009}'s method by exploring the values of the slack variables in order to skip redundant iterations, becoming computationally more efficient. \cite{ZhangReimann2014} addressed the requirement to have the true nadir points, which are often difficult to compute \citep[see][]{AlvesCosta2009,EhrgottRyan2003,KirlikSayin2015}. The proposed methodology skips the redundant iterations in a way that does not require more computations when using approximated nadir values. The disadvantage, however, is that it can only be used for the computation of the whole Pareto front and not for the generation of a representation. 

\cite{EusebioEtAl2014} addresses the representation problem using the $\epsilon$-constraint method for bi-objective problems. \cite{EusebioEtAl2014} propose two algorithms, one for coverage and another for uniformity. The former is an insertion based method: at each iteration, the algorithm adds a solution in between the two most distant in the representation. The latter successively adds solutions spaced by a predetermined step to the representation. Both methods terminate when a desired level of coverage and uniformity are met. 

In this work we take as baseline the work of \cite{MavrotasFlorios2013} and combine it with the strategies presented in \cite{ZhangReimann2014}, which not only allow to pass over redundant iterations, but also permits to overcome the need of computing the true nadir points without requiring the computation of the whole Pareto front. In this way, we can  address the aforementioned drawbacks of both works, while keeping their advantages. Furthermore, we develop three search strategy algorithms, one for coverage, another for uniformity, and a third one for cardinality, addressing each dimension of the Pareto front representation problem. The first two algorithms are based on the work by \cite{EusebioEtAl2014} and are extended for MOILP problems with more than two objectives. The latter refines the virtual grid, introduced by \cite{Mavrotas2009}, after finding each solution, trying to verify the cardinality level by focusing the search on the feasible region. This work, by the very nature of its model formulation and strategy to look for new solutions in the Pareto front, is independent of both the number of non-dominated solutions found and the objective functions ranges. 

The remainder of this paper is organized as follows. Section 2 introduces the mathematical background, namely the type of problem addressed in this paper, the model used and the representation problem. Section 3 presents the proposed methodology, both the generic algorithm and the three proposed search strategies, as well as an illustrative example for each one of them. Section 4 provides the computational results for the strategies presented in Section 3, both the ones for computing the full Pareto front as well as those for the representation problem are presented and discussed. At last, Section 5 presents some concluding remarks and future work lines are put forward.

\section{Mathematical background}
\noindent This sections presents the main concepts, definitions, and notation on multi-objective optimization, the $\epsilon-$constrain approach, and some fundamental concepts related to the representation of the whole set of outcome vectors. 

\subsection{The multi-objective integer programming problem}
\noindent Consider the following multi-objective integer linear programming model. 

\begin{equation}\label{MOO}
    \begin{array}{rcl}\tag{$P$}
        \max z_1(x) & = &  (c^1)^\top x \\
        \vdots      &   &  \vdots \\
        \max z_k(x) & = &  (c^k)^\top x \\
        \vdots      &   &  \vdots \\
        \max z_p(x) & = &  (c^p)^\top x \\
                    &   &  \\
        \mbox{subject to:} &   & x \in X.
    \end{array}
\end{equation}

\noindent where $x = (x_1,\ldots,x_j,\ldots,x_n)$ is an $n-$vector of non-negative and integer \textit{decision variables}, $(c^k)^\top = (c_{1}^k,\ldots,c_{j}^k,\ldots,c_{j}^k)$ is an $n$ row vector composed of the \textit{coefficients} of the decision variables in the \textit{objective functions}, $k=1,\ldots,p$ (we assume these coefficients are integer values or can be converted into integer), and $X$ is the 
\textit{feasible region} in the \textit{decision space}, $\mathbb{Z}^{n}$. Let $Z = z(X)$ denote the image of the \textit{feasible region} according to the objective functions. The set $Z$ is called the \textit{feasible region} in the \textit{objective space}, i.e., $\mathbb{Z}^{k}$ along with the order relation imposed by the objective functions in this set.

Problem \ref{MOO} can be presented in a more compact form as follows. 

\begin{equation}\label{eq:compact}
    \begin{array}{rcl}\tag{$P^{c}$}
        \mbox{``}\max\mbox{''} z(x) & = &  Cx \\
                 \mbox{subject to:} &   & x \in X.
    \end{array}
\end{equation}

\noindent where ``$\max$'' means that all functions are to be maximized, $z(x) = \big(z_1(x),\ldots,z_k(x),\ldots, z_p(x)\big)$ is the vector of the $p$ objective functions, and  $C$ is an $n\times k$ matrix, each line being composed of the coefficients of each objective function.  

\begin{definition}{(Dominance)}
Let $z^\prime$ and $z^{\prime\prime}$ denote two outcome vectors, or solutions, in the objective space. The vector $z^\prime$ \textit{(weakly) dominates} $z^{\prime\prime}$, iff $z^{\prime} \geq z^{\prime\prime}$ and $z^{\prime} \neq z^{\prime\prime}$ (i.e., $z_{k}^{\prime} \geqslant z_{k}^{\prime\prime}$, with at least a strict inequality, for $k=1,\ldots,p$).  When all the inequalities are strict it is said that $z^\prime$ \textit{strictly dominates} $z^{\prime\prime}$.  
\end{definition}

Let $z^{\ast}$ denote a feasible outcome vector,  i.e., $z^{\ast} \in Z$. This vector is a \textit{weak non-dominated outcome vector} if an only if there is no another $z \in Z$ such that $z$ strictly dominates $z^{\ast}$. Some of weak non-dominated outcome vectors are of no interest since we can find another weak non-dominated outcome vector with better performances in at least one objective function (they can thus be discarded). If there is no vector $z \in Z$, such that $z$ weakly dominates $z^{\ast}$, the vector $z^{\ast}$ is simply called a \textit{non-dominated outcome vector}. Whenever $z^\ast$ can be obtained as a weighted-sum of the $p$ objective functions with strictly positive weighing factors, $z^\ast$ is called a \textit{supported outcome vector}. Otherwise, $z^\ast$ is said to be an \textit{unsupported outcome vector}. Let $N(Z)$ denoted the whole set of outcome non-dominated vectors or solutions, also called \textit{Pareto front}. 

The same kind of concepts can be applied in the decision space. Thus, a feasible solution $z^{\ast}$ is called a \textit{weakly efficient solution} if and only if we cannot find another $x \in X$ such that $Cx \geq Cx^\ast$ and $Cx \neq Cx^\ast$. The concepts of \textit{efficient solutions}, \textit{supported efficient solutions}, and \textit{unsupported efficient solutions} are easy to define. Finally, let $E(X)$ denote the set of all efficient solutions. For more details about these definitions see \cite{Ehrgott2005} and \cite{Steuer1986}.


\subsection{The $\epsilon-$constraint approach}
\noindent This subsection is devoted to present an approach for identifying the whole Pareto front. This approach is based on the resolution of a sequence of problems of the following form.

\begin{equation}\label{e-MOO}
    \begin{array}{rl}\tag{$P^\epsilon$}
        {\displaystyle \max_{x}}   & \big\{z_q(x)\big\}   \\
        \mbox{subject to:} &   z_k(x) \geqslant \epsilon_k, \;\; k = 1,\ldots,p,\; k \neq q\\
                      & x \in X.
    \end{array}
\end{equation}

\begin{theorem}{(\citeauthor{HaimesEtAl1971}, 1971)}\label{Theo_HaimesEtAl71}
If $x^\ast$ is an optimal solution of Problem \ref{e-MOO}, for some $q$, then $x^\ast$ is a weakly efficient solution of Problem \ref{MOO}.  
\end{theorem}

\begin{proof}
See  \cite{ChankongHaines2008}, \cite{Ehrgott2005}, or \cite{HaimesEtAl1971},. 
\end{proof}

\vspace{0.5cm}

Algorithm \ref{alg:e-constraint_algo} can be implemented for solving a sequence of \ref{e-MOO} instances. This algorithm requires as input the data or Problem \ref{eq:compact}, $C$ and $X$, and a small enough strictly positive parameter value, $\eta$, and it provides as output the whole Pareto front, $N(Z)$. This algorithm makes use of four internal procedures:

\begin{enumerate}[label={--}]
    \item $Ideal(P)$: which computes the ideal point, $z^\ast = (z^{\ast}_1,\ldots,z^{\ast}_k,\ldots,z^{\ast}_p)$, trough an individual maximization of each objective function. 
    \item $ApproxNadir(P)$: which computes the nadir point or an approximation of it, $z^{nad} = (z^{nad}_1,\ldots,z^{nad}_k,\ldots,z^{nad}_p)$.
    \item From the previous two calculations we can identify the ranges of possible values, for each objective function, which are bounded from below by $z^{nad}_k$ and from above by $z^{\ast}_k$, for $k=1,\ldots,p$.   
    \item $Solve(P^\epsilon)$: which makes use of an integer linear programming solver to solve $P^\epsilon$, and provides the outcome vector $\bar{z} = (\bar{z}_1,\ldots, \bar{z}_k, \ldots, \bar{z}_p)$.
    \item $Filter(\hat{N}(Z))$: which makes the filtering of an auxiliary set, $\hat{N}(Z)$, which may contain weakly non-dominated outcome vectors, and provides as output the Pareto front set, $N(Z)$.
\end{enumerate}

\begin{algorithm}[!htbp]
    \mbox{\bf{Input}:}{ $C,X, \eta$}\;
    \mbox{\bf{Output}:}{ $N(Z)$}\;
    $\hat{N}(Z) \gets \{\}$\;
    $z^\ast \gets Ideal(P)$\;
    $z^{nad} \gets ApproxNadir(P)$\;
    Select $q$ from $\{1,\dots,p\}$ and build $P^\epsilon$\; 
    \For{$(k=1,\; k \neq q)$ \mbox{\bf{to}} $(p)$}
    {
        $\epsilon_k \gets z_{k}^{nad} + \eta$, ~for $k=1,\ldots,p$\;
        $\hat{z}_k \gets -\infty$\;
        \While{$(\hat{z}_k < z_{k}^{\ast})$} 
        {
         	 $\hat{z} \gets Solve(P^\epsilon)$\;
         	 $\hat{N}(Z) \gets \hat{N}(Z) \cup \{\hat{z}\}$\;
         	 $\epsilon_k \gets \hat{z}_k + \eta$\;
         }
    }
    $N(Z) \gets Filter(\hat{N}(Z))$\;
    \mbox{\bf{return}}{$({N}(Z))$}\;
\caption{Computing the Pareto front, $({N}(Z))$.}   
\label{alg:e-constraint_algo}
\end{algorithm}

\vspace{0.5cm}

The $\epsilon-$constraint approach consists of solving a sequence of \ref{e-MOO} instances by adjusting successively (increasing) the parameters $\epsilon_k$, for $k=1,\ldots,p$ with $k \neq q$.  Any basic algorithm designed for implementing this approach (as for example Algorithm \ref{alg:e-constraint_algo}) suffers from three major drawbacks, independently of the solver used for optimizing an instance of \ref{e-MOO}:

\begin{enumerate}
    \item (Model structure) It leads to solve unnecessary instances of Problem \ref{e-MOO}, since it may find some weakly efficient solutions with no interest (i.e.,  disposable weak efficient solutions). This is due to the shape of the model objective function and appears as a conclusion of Theorem \ref{Theo_HaimesEtAl71}.  
    \item (Model structure) It leads to a large amount of CPU time \citep{EhrgottRyan2003}. This is due to the nature of the $\epsilon$ constraints; they lack of some flexibility \citep{EhrgottRyan2003, EhrgottRuzika2008}. 
    \item (Model technical parameter adjustment) It may miss some weakly efficient solutions of interest. This is due to the fact that the solutions obtained are largely dependent on the values chosen for the $\epsilon$ vector \citep{LaumannsEtAl2006}.  
\end{enumerate}

The third drawback was addressed by \cite{LaumannsEtAl2006}. These authors shown in detail the serious issues related with the \textit{a priori} technical adjustment of the $\epsilon$ vector, and propose an adaptive based algorithm with an interesting scheme for determining the adequate values for the $\epsilon$ parameters, when sequentially solving each instance of Problem \ref{MOO}. 

The second drawback was reported and firstly addressed by \cite{EhrgottRyan2003}. These authors proposed an elastic technique for changing the nature of $\epsilon$ type constraints. As a consequence the following model was proposed (some improved versions of this model can be seen in \citealt{EhrgottRuzika2008}). 

\begin{equation}\label{e2s-MOO}
    \begin{array}{rl}\tag{$P^{\epsilon e}$}
         {\displaystyle \max_{x,e^{-},e^{+}}}  & {\displaystyle \left\{ z_q(x)  - \sum_{k=1, \, k\neq q}^{p}p_ke_{k}^{-} \right\}}\\
        \mbox{subject to:} &    z_k(x) + e_{k}^{-} - e_{k}^{+} = \epsilon_k, \;\; k = 1,\ldots,p,\; k \neq q\\
                    &   e_{k}^{-}, e_{k}^{+} \geqslant 0, \;\; k = 1,\ldots,p,\; k \neq q\\
                    &   x \in X.
    \end{array}
\end{equation}

The following result did not change the nature of the output of any algorithm, and the first drawback persists. 

\begin{theorem}{(\cite{EhrgottRyan2003})}
If $p_k >0$, for $k=1,\ldots,p$ with $k \neq p$, and $(x^\ast,e^{-\ast},e^{+\ast})$ is an optimal solution of Problem \ref{e2s-MOO}, for some $q$, then $x^\ast$ is a weakly efficient solution of Problem \ref{MOO}.  
\end{theorem}

\begin{proof}
See \cite{EhrgottRyan2003}, or \cite{EhrgottRuzika2008}
\end{proof}

The first drawback was properly addressed by \cite{Mavrotas2009} with the proposal of a new model. This new model is based on the introduction of slack variables and it seems to mitigate also the second drawback. A new improvement for multi-objective integer linear programming which makes a slight change on the objective function was proposed by \cite{MavrotasFlorios2013}. This change will not alter the theoretical results presented in \cite{Mavrotas2009}. The improved model can be stated follows. 

\begin{equation}\label{es-MOO}
    \begin{array}{rl}\tag{$P^{\epsilon s}$}
        {\displaystyle \max_{x,s}}  & {\displaystyle \left\{z_q(x)  + \rho\sum_{k=1, \, k\neq q}^{p}10^{k-1}\frac{s_k}{r_k}\right\}}\\
        \mbox{subject to:} &    z_k(x) - s_{k}  = \epsilon_k, \;\; k = 1,\ldots,p,\; k \neq q\\
                    &   s_{k},  \geqslant 0, \;\; k = 1,\ldots,p,\; k \neq q\\
                    &   x \in X.
    \end{array}
\end{equation}
\noindent where $r_k$ is the amplitude of the range values for each objective function, $k = 1,\ldots,p,\, k \neq q$.

The nature of the output of any algorithm designed to solve the previous problem changes as it can be see in the following theoretical result.

\begin{theorem}{(\citealt{Mavrotas2009})}
If $(x^\ast,s^\ast)$ is an optimal solution of Problem \ref{es-MOO}, for some $q$ and for $\rho >>0$, then $x^\ast$ is an efficient solution of Problem \ref{MOO}.  
\end{theorem}

\begin{proof}
See \cite{Mavrotas2009}.
\end{proof}

Problem \ref{es-MOO} will be used in our algorithms. 

\subsection{On the representation of the Pareto front}\label{Background representation}
\noindent A discrete representation of the Pareto front, $R(N)$, is a finite subset of the Pareto front, $N(Z)$. The quality of representations, in terms of how well its subset captures the characteristics of the full set, has been the focus of many studies \citep[see][]{Sayin2000,FaulkenbergWiecek2010,AudetEtAl2021}. These studies proposed dimensions of interest in their evaluation of the representation. \cite{Sayin2000} proposes that the quality of solution representations should be assessed through three criteria: 
\begin{enumerate}
\item \textit{Coverage}: How well does the representation covers all regions of the objective space $Z$ included in $N(Z)$.
\item \textit{Uniformity}: How diverse and equally spaced are the points in the representation, i.e., how spread are the points (outcome vectors) included in the representation.
\item \textit{Cardinality}: The number of outcome vectors considered in the representation, $\Pi\big(R(N)\big)$. 
\end{enumerate}

Coverage and uniformity are dependent on the distance between points in the representation. There are multiple distance metrics, the most common being as follows:

\begin{equation*}
    d(z,z') = \begin{cases}
    \left( |z_1-z'_1|^t +\dots+|z_k-z'_k|^t +\dots +  |z_p-z'_p|^t\right)^{1/t}, & t \geqslant 1, \\
    ~~ \\
    \max \left(|z_1-z'_1|,\dots,|z_k-z'_k|,\dots, |z_p-z'_p|\right), & t = \infty.
    \end{cases}
\end{equation*}

\noindent where, if $t=1$ the distance metric used is the Manhattan distance, if $t=2$ it corresponds to the Euclidean distance, and if $t=\infty$ the Chebyshev distance is applied. The higher the $t$ value, the smaller the distance value measured. Considering a distance metric, the coverage error of the representation can be stated as follows:
\begin{equation}
    \Gamma\big(R(N), N(Z)\big) = \max_{z\in N(Z)} \min_{z'\in R(N)} d(z,z'),
    \label{eq: coverage_gamma_gen}
\end{equation}
which corresponds to the maximum distance from $z\in N(Z)$ to its closest point in the representation $z'\in R(N)$. For the sake of simplicity, from hereinafter, let us use the notation $\Gamma\big(R(N)\big)$ instead of $\Gamma\big(R(N), N(Z)\big)$. A representation $R(N)$ is said to have coverage $\gamma$, if $\Gamma\big(R(N)\big)\leqslant\gamma$. To increase coverage, the coverage error $\Gamma\big(R(N)\big)$ must be minimized. Since the representation $R(N)$ is a finite subset of $N(Z)$, minimizing the maximum distance from $z\in N(Z)$ to its closest point in the representation $z'\in R(N)$ can be done by minimizing the maximum distance between each two consecutive points in the representation. Furthermore, by guaranteeing that all consecutive points in the representation $R(N)$ are distanced by a value lower or equal to $\gamma$, there is also the guarantee that $R(N)$ has, at most, a coverage error of $\gamma$.

The uniformity level corresponds to the minimum distance between any two points in the representation:
\begin{equation}
    \Delta\big(R(N)\big) = \min_{z,z'\in R(N), z\neq z'}d(z,z')
\end{equation}

A representation $R(N)$ is said to have uniformity $\delta$, if $\Delta\big(R(N)\big)\geqslant \delta$. To increase uniformity, the uniformity level $\Delta\big(R(N)\big)$ should be maximized. As a consequence of the aforementioned, the discrete representation problem, \ref{DRP}, can be stated as a three-objective optimization problem \citep{ShaoEhrgott2016}:

\begin{equation}\label{DRP}
    \begin{array}{l}\tag{$P^{DRP}$}
        \min \; \Pi\big(R(N)\big) \\
        \max \; \Delta\big(R(N)\big) \\
        \min \; \Gamma\big(R(N)\big)  \\
        \mbox{subject to:} \;\; R(N) \subset N(Z), |R(N)|>\infty 
    \end{array}
\end{equation}

The objectives considered in Problem \ref{DRP} are interrelated. \cite{Sayin2000} noted that, by improving the representation coverage, cardinality would increase too. As uniformity increases, potentially, the opposite effect on cardinality may emerge. Furthermore, \cite{KiddEtAl2020} proved that, when comparing two representation, $R^{\prime}(N)$ and $R^{\prime\prime}(Z)$, with the same cardinality, where $R^{\prime\prime}(Z)$ is an equidistant representation, then $\Gamma\big(R^{\prime\prime}(Z)\big)\leqslant\Gamma(R^{\prime}(N))$ and $\Delta\big(R^{\prime\prime}(Z)\big)\geqslant\Delta(R^{\prime}(N))$. 

In this work, we make use of the Chebyshev norm to compute $\Gamma(R(N))$ and $\Delta(R(N))$. Although other metrics may be used, the choice is supported by three main aspects: (1) Chebyshev norm allows the decoupling of the distance for each criterion, in the sense that all coordinate distances will be at most the Chebyshev distance, i.e., if $|z_k - z'_k|\leqslant 1$  for all $k = 1,\ldots, p$, then $d_{l^\infty}(z,z')\leqslant 1$; (2) From the previous fact, it also derives that both the coverage error and the uniformity level have a direct correspondence to the coordinate distances without requiring the combination of multiple coordinates; this allows greater control over the representation criteria; and (3) The combination of coordinates, a fundamental characteristic of the other metrics, may distort the distance metric, consequently the coverage error and uniformity level, if the ranges of each criterion are significantly different \citep{Sayin2000}.

\section{Representation of the Pareto Front}
\begin{figure}[!ht]
    \centering
    \includegraphics[width=0.95\textwidth]{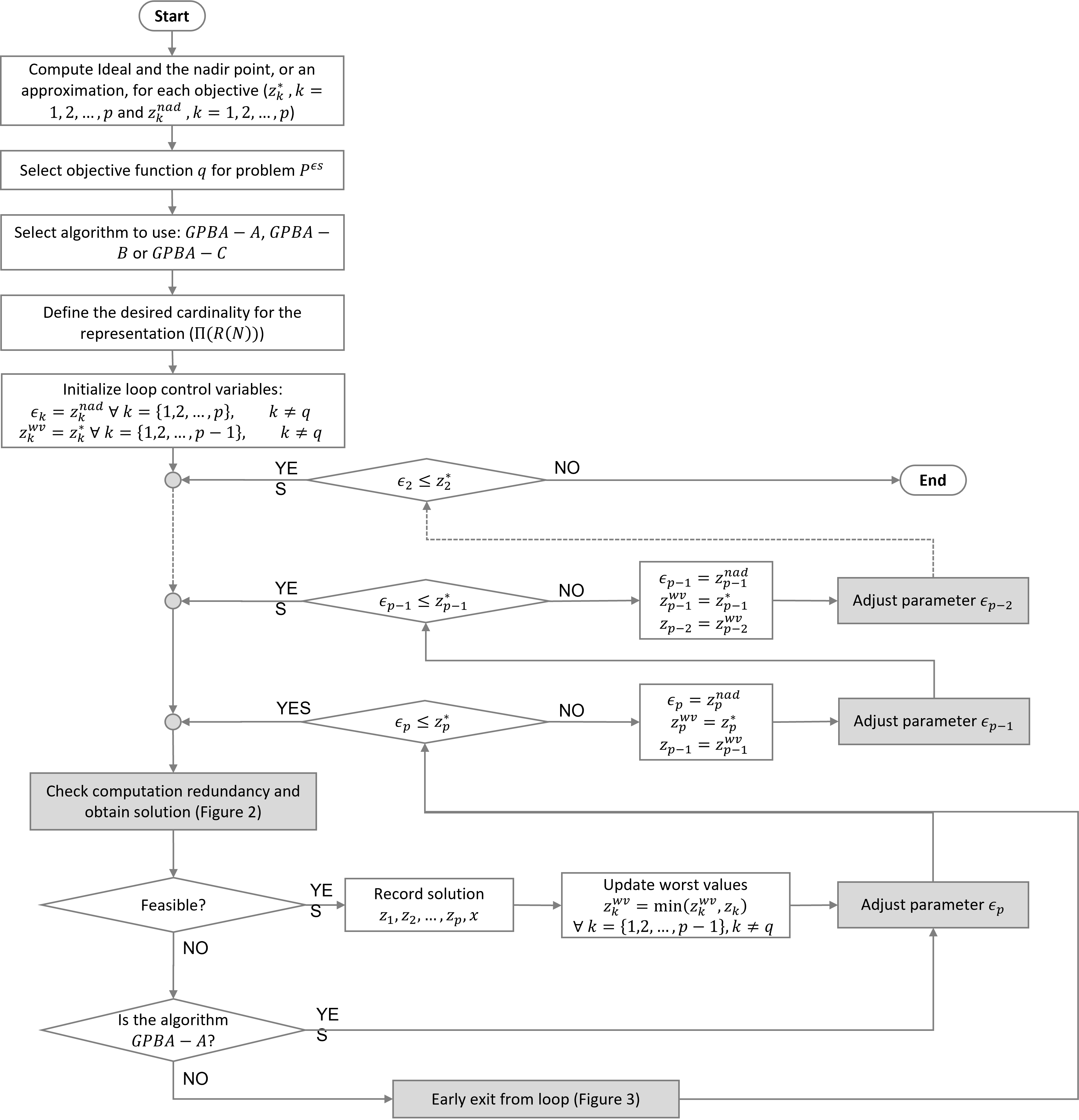}
    \caption{Flowchart for the generation algorithm}
    \label{fig:flow_generic}
\end{figure}

\noindent This section presents the three algorithms developed for finding a representation of the Pareto front for MOILP Problems \ref{MOO}: \textit{GPBA-A}, \textit{GPBA-B} and \textit{GPBA-C}. Each algorithm targets one dimension of the discrete representation Problem \ref{DRP} (coverage, uniformity, and cardinality) by employing different strategies when exploring the feasible region. The first two algorithms, \textit{GPBA-A} and \textit{GPBA-B}, are an extension of the ones proposed by \cite{EusebioEtAl2014}, aiming to improve coverage and uniformity, respectively. The third algorithm, \textit{GPBA-C}, represents an improvement over \cite{MavrotasFlorios2013} algorithm, and concerns mostly cardinality. However, all algorithms share the same main procedure, detailed under Algorithm \ref{alg:e-constraint_algo} and are represented in Figure \ref{fig:flow_generic}, whilst solving Problem \ref{es-MOO}.

\begin{figure}[!ht]
    \centering
    \includegraphics[width=0.95\textwidth]{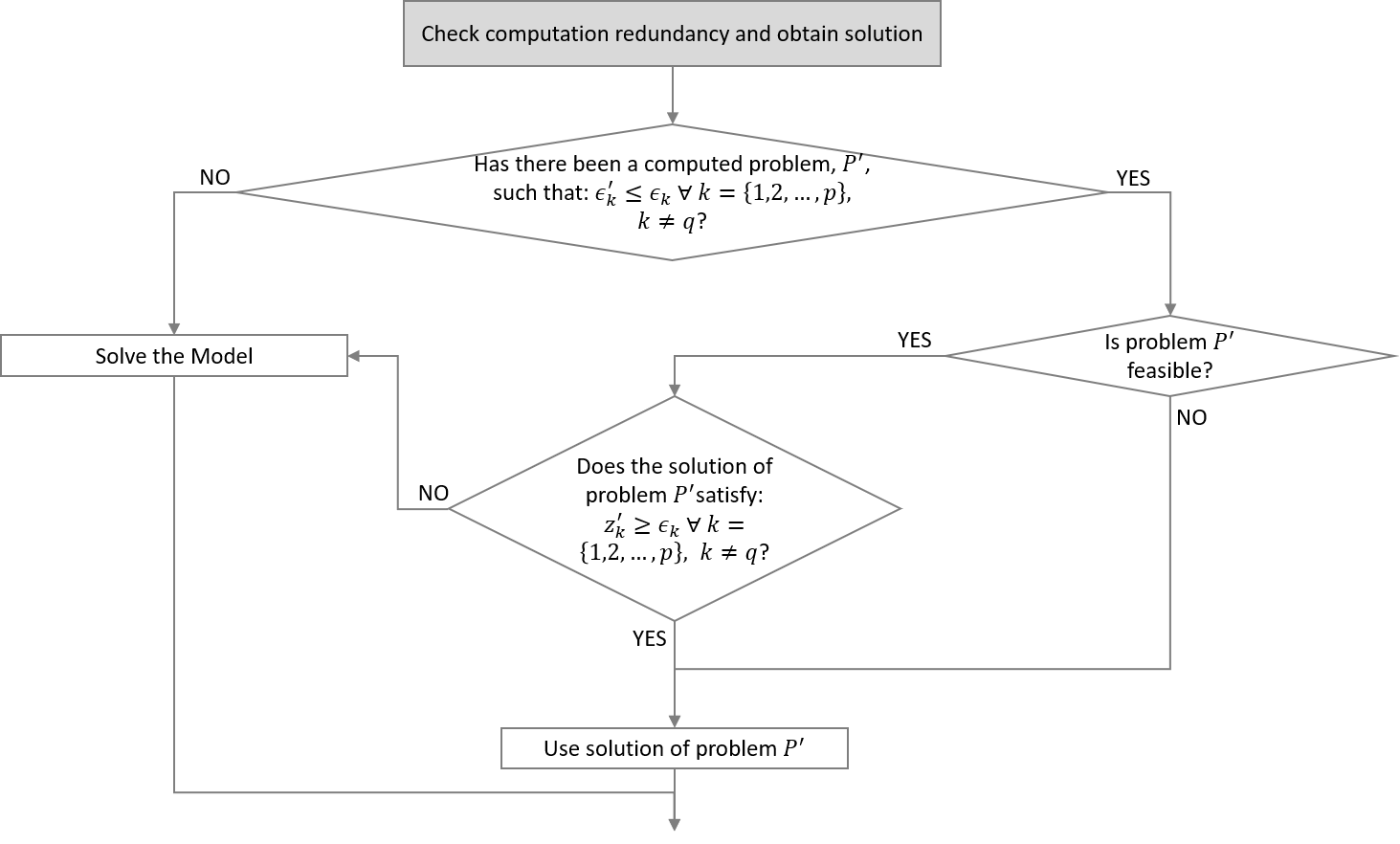}
    \caption{Process to determine if the problem needs to be solved or can be skipped.}
    \label{fig:flow_check_redundancy}
\end{figure}

\begin{figure}[!ht]
    \centering
    \includegraphics[width=0.8\textwidth]{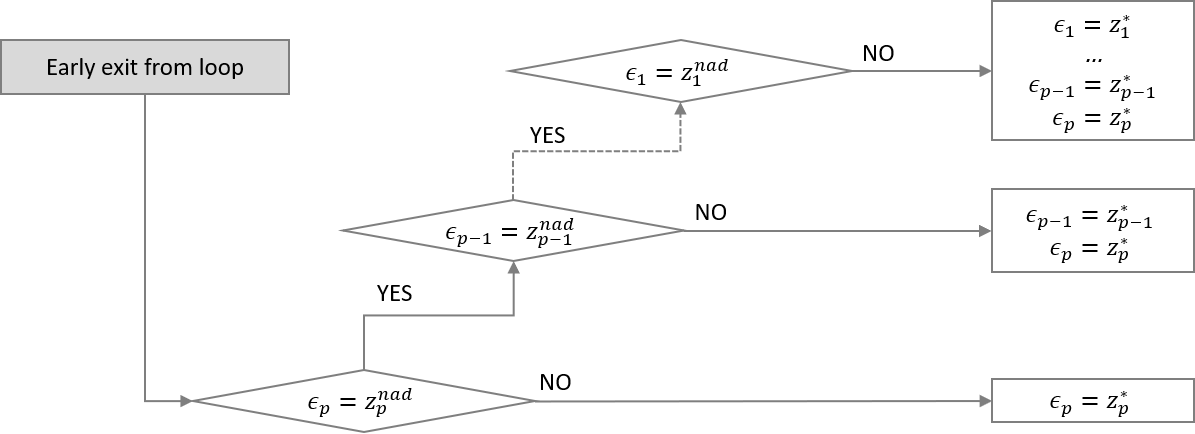}
    \caption{Accelerated exit algorithm.}
    \label{fig:early_exit}
\end{figure}

\begin{enumerate}
    \item The first step of the algorithm consist of the computation of the lower (pessimistic) and upper (optimistic) bounds of the Pareto front, $z^{nad}$ a and $z^\ast$, respectively. Although the optimistic value is typically easily computed by maximizing  each objective function separately, the inverse cannot always be done to obtain the pessimistic values. Many studies in the literature show the difficulty of obtaining good estimations for those bounds \citep[see][]{IsermannSteuer1988,AlvesCosta2009}. Consequently, nadir point overestimation is common. All three proposed methods are able to handle overestimated bounds, allowing any strategy to compute the pessimistic values for each objective function. 
    \item The objective function $q$ to be directly considered in Problem \ref{es-MOO} is then chosen, as well as the representation algorithm (\textit{GPBA-A}, \textit{GPBA-B} or \textit{GPBA-C}). Depending on the representation algorithm chosen, the desired characteristics of the representation must be defined: coverage error, uniformity level or the maximum cardinality of the representation.
    \item The forth step consists of the initialization of the main loop by enforcing the parameters $\epsilon_k$, in Problem \ref{es-MOO}, to be the same as the pessimistic values for each objective with a small perturbation, resulting in \ref{es-MOO} being in its most relaxed form. Additionally, the relative worst values for each objective function, $z^{wv}_k$, are also initialized as the ideal values.
    \item The third step, consists of the comparison of the problem's parameters against the results of the previously solved problem, in order to check if a new outcome vector needs to be computed or if it is going to yield a redundant result. That procedure is shown in Figure \ref{fig:flow_check_redundancy}:
    \begin{enumerate}
        \item Select from the list of previous iterations, for each $k=1,\ldots,p,\; k\ne q$, a Problem $P'$ such that the $\epsilon'$ vector is equal to $\epsilon$ or closer to the the nadir point, i.e., $z^{nad}_k\leqslant \epsilon'_k \leqslant \epsilon_k,$ for all $k=1,\ldots,p,\; k\ne q$.
        \item If such a problem does not exist or, if it exists, its solution does not lay in the bounds of Problem \ref{es-MOO} when using $\epsilon$, then solve Problem \ref{es-MOO} and return to the main procedure (see Figure \ref{fig:flow_generic}).
        \item If such a Problem $P'$ exists but it is infeasible, then return to the main procedure (see Figure \ref{fig:flow_generic}), considering the problem as an infeasible one.
        \item If such a Problem $P'$ exists and it is feasible, and its solution is within the bounds of Problem \ref{es-MOO}, when using $\epsilon$, i.e., $z'_k \geqslant \epsilon_k$, for all $k=1,\ldots,p,\; k\ne q$, then return to the main procedure (see Figure \ref{fig:flow_generic}) considering the solution of Problem $P'$ as the solution of Problem \ref{es-MOO}.
    \end{enumerate}
    \item If the obtained (or selected from the set of previously computed solution) outcome vector is feasible, such a vector is saved and the new $z^{wv}_k$ updated for all objective functions, with the exception of function $q$ and the innermost loop $p$. Then, the new $\epsilon_p$ parameter is readjusted according to the algorithm used. Otherwise, and in case the representation algorithm is not \textit{GPBA-A}, the early exit loop is used (see Figure \ref{fig:early_exit}).
    \item The process is repeated until the optimistic values is reached in all loops, i.e., until Problem \ref{es-MOO} reaches its most constrained form.
\end{enumerate}

The aforementioned acceleration strategies, the redundancy checking (Figure \ref{fig:flow_check_redundancy}) and the early exit from loop (Figure \ref{fig:early_exit}) procedures, were originally proposed by \cite{ZhangReimann2014}. While the latter is only based on Proposition \ref{prop: infeasibility}, the former resorts to both Proposition \ref{prop: infeasibility} and Proposition \ref{prop: optimality}. Let $LP^{\epsilon s}$ denote a linear relaxation of Problem \ref{es-MOO}:

\begin{proposition}{(Infeasibility)}\label{prop: infeasibility}
If Problem $LP^{\epsilon s}$ is infeasible, then Problem \ref{es-MOO} is also infeasible.
\end{proposition}

\begin{proposition}{(Optimality)}\label{prop: optimality}
If $x^*$ is the optimal solution of Problem $LP^{\epsilon s}$, and $x^*$ is in the feasible region of \ref{es-MOO}, then $x^*$ is also the optimal solution of Problem \ref{es-MOO}.
\end{proposition}

In the following subsections, we will details the algorithms to compute the parameters $\epsilon_k$, for Problem \eqref{es-MOO}.

\subsection{Coverage representation  (\textit{GPBA-A})}
\noindent When looking for a representation that privileges the coverage of the Pareto front, the algorithm will aim at guaranteeing that each area of the Pareto front is represented. This may be achieved, as discussed in Section \ref{Background representation}, by minimizing the maximum distance between any two consecutive points in the representation ($\Gamma\big(R(N)\big)$), Equation \ref{eq: coverage_gamma_gen}, i.e. the coverage error. The quality of the representation is controlled by the parameter $\gamma$, such that $\Gamma\big(R(N)\big)\leqslant \gamma$.

Parameter $\gamma$, corresponding to the acceptable coverage error, is sometimes difficult to determine. An alternative to directly conveying the acceptable coverage error, is to define an acceptable cardinality in each objective ($c_k$, for all $k = {1,\dots,p}$, $k\neq q$) based on their range ($z^*_k - z^{nad}_k$, for all $k = {1,\dots,p}$, $k\neq q$). The ratio between the range and cardinality provides the acceptable coverage error for each objective: $\gamma_k = (z^*_k - z^{nad}_k)/c_k$. The sum of the specified cardinality in each objective will correspond to the approximate cardinality of the representation, i.e., $\Pi\big(R(N)\big) = \sum_{k=1,\, k\neq q}^{p} c_k$. Note that, the cardinality of the representation $\Pi\big(R(N)\big)$ is merely an approximation since the distribution of the Pareto front is unknown and it is not this parameter that controls the representation.

Algorithm \ref{alg:GPBA-A} presents the method to adjust parameter $\epsilon_k$ when aiming to have a representation that privileges coverage. It performs a single run for every iteration and for each objective $k$, i.e. for each loop, after solving a problem and obtaining, or not, a new outcome vector. The algorithm takes six inputs: (1) the desired coverage error over that loop, $\gamma_k$; (2) the value of the $\epsilon_k$ parameter of Problem \ref{es-MOO} used in the current iteration; (3) the $k^{th}$ component of the outcome vector $z$, $z_k$, for the Problem \ref{es-MOO} used in the current iteration; (4) the $k^{th}$ component of the ideal vector $z^{\ast}$, $z^{\ast}_k$; (5) the $k^{th}$ component of the approximation to the nadir vector $z^{nad}$, $z^{nad}_k$; and (6) the set $D$, corresponding to the ordered set of the discarded points for the current loop, which contains the points already obtained as well as the areas where redundant points will be obtained. The algorithm's outputs are the updated value of parameter $\epsilon_k$ to be used in the next iteration, and the updated ordered set of discarded points, $D$, for the current loop. To this end, the first part of the algorithm updates $D$ based on one of the following two cases:

\begin{enumerate}
    \item If the Problem \ref{es-MOO} solved in this iteration, when using $\epsilon_k$, is infeasible (the feasible region is empty, $X^{^\epsilon s} = \{\}$), then if the problem is more constrained it will be infeasible too, as stated in Proposition \ref{prop: infeasibility}. As a result, all the points ranging from $\epsilon_k$ to the $k^{th}$ component of the ideal vector, $z^{\ast}_k$, can be discarded, and should be added to $D$ for the current loop.
    \item If the Problem \ref{es-MOO} solved in this iteration, when using $\epsilon_k$, is feasible, all the points ranging between $\epsilon_k$ to the $k^{th}$ component of the obtained vector, $z_k$, can be discarded and should be added to $D$. This is based on Proposition \ref{prop: optimality}. 
\end{enumerate}

Afterwards, using the updated set of discarded points, $D$, the algorithm computes the new value for parameter $\epsilon_k$:

\begin{enumerate}
    \item In the case where the current value of $\epsilon_k$ is the $k^{th}$ component of the pessimistic vector, $z^{nad}_k$, the problem solved in this iteration corresponds to the first extreme point. Hence, the following point to compute is the next extreme, the ideal value $z^{\ast}_k$.
    \item Otherwise, $D$ will have at least more than two points. As a result, the two consecutive points distancing the most from each other are drawn from set $D$. Then, on the one hand, if the distance between those points meets the desired coverage error, the loop may be exited. To do so, set $D$ is updated as empty and the $\epsilon_k$ parameter is set to a value higher than the ideal $z^{\ast}_k + 1$. On the other hand, if the distance between those points meets the desired coverage error, $\epsilon_k$ is set to find a point that minimizes the coverage error, i.e., a point at least in between the two most distant consecutive points, $(z_1 - z_2)/2$.
\end{enumerate}

\begin{algorithm}[!ht]
    \mbox{\bf{Input}:}{ $\gamma_k$, $\epsilon_k$, $z_k$, $z^{\ast}_k$, $z^{nad}_k$, $D$}\;
    \mbox{\bf{Output}:}{ $\epsilon_k$, $D$}\;
    \eIf{$(X^{\epsilon s} = \{\})$}{
    $D \gets D \cup \{\lceil\epsilon_k\rceil,\lceil\epsilon_k\rceil+1,\dots,z^*_k\}$\;}
    {
    $D \gets D \cup \{\lceil\epsilon_k\rceil,\lceil\epsilon_k\rceil+1,\dots,z_k\}$\;}
    \eIf{$(\epsilon_k = z^{nad}_k)$}
    {$\epsilon_k \gets z^*_k$\;}{
    {Let $z_1$ and $z_2$ be the most distant consecutive values in $D$}\;
     \eIf{$(d(z_1,z_2) \leqslant \gamma_k)$}{$\epsilon_k \gets z^*_k+1$\;
     $D \gets \{\}$\;}
    {$\epsilon_k \gets (z_{1} + z_{2})/2$\;}}
    \mbox{\bf{return}}{$(\epsilon_k$, $D)$}\;
\caption{Adjusting the parameter $\epsilon_k$ in \textit{GPBA-A}.}   
\label{alg:GPBA-A}
\end{algorithm}

\subsection{Uniformity representation (\textit{GPBA-B})}
\noindent A representation that privileges uniformity is a representation which spread points as equally spaced as possible. By maximizing the minimum distance between two points in the representation, the uniformity level $\Delta\big(R(N)\big)$ is increased. This representation is controlled by an acceptable uniformity level $\delta$.

Similar to the acceptable coverage level, the acceptable uniformity level may be difficult to define. Hence, it is advised to derive it based on the range of each objective function values, ($z^*_k - z^{nad}_k$, for all $k = 1,\dots,p$, $k\neq q$) and the acceptable cardinality in each objective ($c_k$, for all $k = 1,\dots,p$, $k\neq q$). The ratio between the range and cardinality provides the acceptable uniformity level for each objective: $\delta_k = (z^*_k - z^{nad}_k)/c_k$.

The procedure to determine $\epsilon_k$, for the next iteration, is presented in Algorithm \ref{alg:GPBA-B}. It takes as input the defined uniformity level $\delta_k$ and the $k^{th}$ component of the obtained vector, $z_k$. The algorithm in each iteration adds $\delta_k$ to $z_k$, until $\epsilon_k$ is greater than the maximum value for objective $k$, $z^{\ast}_k$. In that case the loop is exited since the problem becomes infeasible. 

\begin{algorithm}[!ht]
    \mbox{\bf{Input}:}{ $\delta_k$, $\epsilon_k$}\;
    \mbox{\bf{Output}:}{ $\epsilon_k$}\;
    {$\epsilon_k \gets z_k + \delta_k$\;}
    \mbox{\bf{return}}{$(\epsilon_k)$}\;
\caption{Adjusting the parameter $\epsilon_k$ in \textit{GPBA-B}.}   
\label{alg:GPBA-B}
\end{algorithm}

\subsection{Cardinality representation (\textit{GPBA-C})}

\noindent The cardinality representation here proposed aims to provide a range insensitive search strategy that seeks to distribute the desired number of solutions in the feasible region. When the desired cardinality is specified, it is very common that the ranges of the objectives in which it is based are overestimated. As a result, the final cardinality of the representation will be significantly lower than the originally considered, and in which the determination of parameters, such as uniformity level $\delta$ and coverage error $\gamma$, may have been based. The cardinality representation algorithm addresses this drawback by defining the search region in a way that intends to maintain the originally determined cardinality, and provides a balance between uniformity and coverage. 

The cardinality representation algorithm starts by defining a uniform grid between the ideal and nadir points. Whenever a solution is computed that skips a step on the grid, the grid is redefined from that point until the ideal value, considering only the remaining points left to compute. This results in a refinement of the search grid within the feasible region of the problem, and consequently in a range insensitive search strategy. The lower the desired cardinality is, the least amount of steps may need to be skipped, hence the more uniform the final representation will be. On the other hand, the higher the desired cardinality, the more solutions may need to be skipped and, consequently, the grid will be refined more times, leading to a representation with a better coverage and worst uniformity. 

Algorithm \ref{alg:GPBA-C} presents the procedure to adapt parameter $\epsilon_k$. It takes six inputs: (1) the grid starting point for that loop, $z^{start}_k$; (2) the $k^{th}$ component of the ideal vector $z^{\ast}$, i.e. the grid end point for that loop, $z^{\ast}_k$; (3) the number of solutions to obtain from that grid, $c'_k$; (4) the position within the grid for this objective, i.e. the grid point number, $i_k$; (5) the $k^{th}$ component of the outcome vector $z$, $z_k$, for the Problem \ref{es-MOO} used in the current iteration; and (6) the $k^{th}$ slack variable, $s_k$, for the Problem \ref{es-MOO}. For the the first iteration, prior to entering Algorithm \ref{alg:GPBA-C}, the grid start point, $z^{start}_k$, corresponds to the nadir approximation vector, $z^{nad}_k$, the position within the grid, $i_k$, to zero, and  the number of solutions to obtain, $c'_k$, to the cardinality for that objective minus 1, $c_k-1$. The algorithm outputs the new parameter $\epsilon_k$ and the updated grid parameters, $z^{start}_k$, $c'_k$ and $i_k$. To that end, the following procedure is applied:
\begin{enumerate}
    \item The first step determines if the outcome vector from the problem solved in the previous iteration, $z_k$, makes any of the next grid points redundant according to Proposition \ref{prop: optimality}, i.e. if a step on the grid may be skipped. Accordingly, the grid step size, $step$, is computed by dividing the grid range by the number of grid points in that grid. The integer part of the division between the slack variable, $s_k$, and the grid step size, $step$, determines how many grid points may be skipped.
    \begin{enumerate}
        \item In case there are steps on the grid to skip, the grid is redefined. The new starting point for the grid, $z^{start}_k$, will become the corresponding component of the outcome vector obtained in the current iteration, $z_k$. The updated number of grid points for the new grid, $c'_k$, corresponds to the number of points from the current grid subtracted the number of grid points already computed, $i_k$. The updated the position within the grid, $i_k$, will be one, i.e. at the beginning of the new grid. At last, the step according to the newly defined grid is computed.
        \item Otherwise, the grid is not updated and the position within the grid, $i_k$, is incremented.
    \end{enumerate}
    \item The second step, makes use of the grid parameters computed in the previous step, $z^{start}_k$, $i_k$ and $step$, the updated parameter $e_k$ is computed. In case $e_k$ surpasses the ideal vector value $z^{\ast}_k$, the loop should be exited and, therefore, the grid parameter should return to their initial values: the grid start point, $z^{start}_k$, corresponds to the nadir approximation vector, $z^{nad}_k$; the number of solutions to obtain, $c'_k$, to the cardinality for that objective, $c_k$; and the position within the grid, $i_k$, to zero.
\end{enumerate}

\begin{algorithm}[!ht]
    \mbox{\bf{Input}:}{ $z^{start}_k$, $z^{\ast}_k$, $c'_k$, $i_k$, $z_k$, $s_k$}\;
    \mbox{\bf{Output}:}{ $\epsilon_k$, $i_k$, $c'_k$, $z^{start}_k$}\;
    {$step \gets \max\big\{(z^{\ast}_k-z^{start}_k)/c'_k,\,1\big\}$\;
    $b \gets \lfloor|s_k/step|\rfloor$\;}
    \eIf{$(b>0)$}{
    $z^{start}_k \gets z_k$\;
    $c'_k \gets c'_k - i_k$\;
    $i_k \gets 1$\;
    $step \gets \max\big\{(z^{\ast}_k-z^{start}_k)/c'_k,\,1\big\}$\;}
    {$i_k \gets i_k +1$\;}
    {$\epsilon_k \gets z^{start}_k + i_k\times step$\;}
    \If{$(\epsilon_k>z^{\ast}_k)$}{
    $z^{start}_k \gets z^{nad}_k$\;
    $c'_k \gets c_k$\;
    $i_k \gets 0$\;}
    
    \mbox{\bf{return}}{$(\epsilon_k$, $z^{start}_k$, $i_k$, $c'_k)$}\;
\caption{Adjusting the parameter $\epsilon_k$ in \textit{GPBA-C}.}   
\label{alg:GPBA-C}
\end{algorithm}

\subsection{An illustrative example}\label{sec: An illustrative example}
\noindent Consider the following numerical example, originally presented in \cite{IsermannSteuer1988}:
\[
    \begin{array}{rcrrcrcrcrcrcrcrcl}
       \max \, z_1(x)   & = &  & 2x_1 &   &  &  & & - & 2x_4 & & & - & 2x_6 & - & 2x_7 & &  \\
       \max \, z_2(x)   & = & - & 2x_1 & + & x_2 & + &  2x_3 & - & x_4 & + & x_5 & + & 2x_6  & - & x_7 & & \\
       \max \, z_3(x)   & = & - & x_1 & - & 2x_2 &  & & - & 2x_4 & + & 3x_5 & + & x_6 & & & &   \\
    \mbox{subject to:}  & & & x_1 & + & x_2 & + & 3x_3 & & & + & 3x_5 & + & 2x_6 & & & \leqslant & 61 \\
     & & & & & 3x_2 & + & 2x_3 & + & 4x_4 & & & & & & & \leqslant & 72 \\
     & & & 5x_1 & + & 3x_2 & & & & & + & 5x_5 & + & 4x_6 & + & 4x_7 & \leqslant & 76 \\
     & & & 4x_1 & + & 2x_2 & & & + & 4x_4 & & & + & 4x_6 & & & \leqslant & 51  \\
     & & & 5x_1 & + & 2x_2 & & & + & 3x_4 & + & x_5 & + & 4x_6 & & & \leqslant & 66 \\
     & & & 2x_1 & + & 2x_2 & & & + & 4x_4 & + & 4x_5 & + & 4x_6 & + & 5 x_7 & \leqslant & 59 \\
     & & & 3x_1 & & & + & 2x_3 & & & + & 5x_5 & + & x_6 & + & 2x_7 & \leqslant & 77 \\ 
     & & & x_1, & & x_2, & & x_3, & & x_4, & & x_5, & & x_6, & & x_7 & \in & \mathbb{Z}_{0}^{+}.\\
    \end{array}
\]

The three representation algorithms, \textit{GPBA-A}, \textit{GPBA-B}, and \textit{GPBA-C}, will be illustrated in this section. To that end, the first iteration over the outermost loop will be used. Common to all representations is the need to compute both the ideal point and the nadir point, or, alternatively, an approximation of the nadir point. In this case, the ideal point is $z^{\ast} = (24, 49, 42)$, and the nadir approximation point, obtained through the minimization of each individual objective function, is $z^{nad} = (-28, -28, -48)$. Additionally, $z_1$ will be used as $z_q$, $z_2$ as the outermost loop, and $z_3$ as the innermost loop. Hence, all algorithms start with $\epsilon_2 = -28$ and $\epsilon_3 = -48$, as well as with the relative worst value for the outermost loop $z^{wv}_2 = 49$.

\subsubsection{Coverage representation (\textit{GPBA-A})}
\noindent Assume the acceptable coverage error for the representation being computed is $\gamma = 15$, for all objective functions. The first iteration over the outermost loop will have $\epsilon_2 = -28$. Then, the iterations over the innermost loop, $k=3$, are the following:
\begin{enumerate}
    \item The first iteration over the innermost loop (Figure \ref{fig:GPBA-A-1}) solves Problem \ref{es-MOO} using $\epsilon_2 = -28$ and $\epsilon_3 = -48$. Point $z^{1}=(24,9,-14)$ is obtained and added to the representation $R(N)=\{z^{1}\}$. Next, the relative worst value for objective function $z_2$, $z^{wv}_2$, and the ordered set of discarded points, $D$, are updated to $z^{wv}_2=9$ and $D=\{-48,-47,\dots,-14\}$. Since the current value for $\epsilon_3$ is $z^{nad}_3$, the new value for $\epsilon_3$ is $42$, corresponding to the ideal value.
    \item Since, according to Figure \ref{fig:flow_check_redundancy}, Problem \ref{es-MOO} needs to be solved for $\epsilon_2 = -28$ and $\epsilon_3 = 42$, the point obtained is $z^{2}=(0,20,42)$ (Figure \ref{fig:GPBA-A-2}). Point $z^{2}$ is added to the representation, resulting in $R(N)=\{z^{1},z^{2}\}$. Since $z^{2}_2 \geqslant z^{wv}_2$, the worst value for objective function $z_2$ does not need to be updated. To the set of discarded points, $z^{2}_3$ is added $D=\{-48,-47,\dots,-14,42\}$. The two most distant consecutive values in $D$ are $-14$ and $42$, hence the new value of $\epsilon_3$ is $(42+(-14))/2 = 14$.
    \item Again, according to Figure \ref{fig:flow_check_redundancy}, Problem \ref{es-MOO} is solved for $\epsilon_2 = -28$ and $\epsilon_3 = 14$, obtaining $z^{3}=(14,13,14)$ (Figure \ref{fig:GPBA-A-3}), and updating the representation to $R(N)=\{z^{1},z^{3}, z^{2}\}$. Again, the worst value for objective function $z_2$, $z^{wv}_2$, does not require updating, while the set of discarded points is updated to $D=\{-48,-47,\dots,-14,14,42\}$. At this point the maximum distance between points $D$ is $28$. As a result, the new value for $\epsilon_3$ is $(14+(-14))/2 = 0$.
    \item Problem \ref{es-MOO} is solved for $\epsilon_2 = -28$ and $\epsilon_3 = 0$, obtaining $z^{4}=(22,6,1)$ (Figure \ref{fig:GPBA-A-4}), and updating the representation to $R(N)=\{z^{1},z^{4}, z^{3}, z^{2}\}$. The worst value for objective function $z_2$, $z^{wv}_2$, needs to be updated $z^{wv}_2 = 6$. The set of discarded points is updated to $D=\{-48,-47,\dots,-14,0,1,14,42\}$. The new value for $\epsilon_3$ is $(42+14)/2 = 28$.
    \item Using $\epsilon_2 = -28$ and $\epsilon_3 = 28$, Problem \ref{es-MOO} is solved and $z^{5}=(8,13,29)$ is obtained (Figure \ref{fig:GPBA-A-5}). Outcome vector $z^{5}$ is then added to the representation, $R(N)=\{z^{1},z^{4}, z^{3}, z^{5},$ $z^{2}\}$. The worst value for objective function $z_2$ remains the same, $z^{wv}_2 = 6$. The set of discarded points is updated to $D=\{-48,-47,\dots,-14,0,1,14,28,29,42\}$. At this point, the maximum distance between points in the set of discarded points is $14< \gamma$. Hence, the loop may stop, and the coverage error obtained for the first iteration over the outermost loop was $13 < \gamma$.
\end{enumerate}

For the next iteration over the loop of objective function $z_2$, the value of $z^{wv}_2 = 6$ is used as $z_2$ for the discarded vector of that loop.

\subsubsection{Uniformity representation (\textit{GPBA-B})}
\noindent Suppose the objective is to have a representation with at least a uniformity level of $\delta = 10$ over all objective functions. The algorithm begins with $e_2=-28$, for the outermost loop, and the following iterations occur over the innermost loop:
\begin{enumerate}
    \item In the first iteration, Figure \ref{fig:GPBA-B-1}, Problem \ref{es-MOO} is solved using $\epsilon_3=-48$. Point $z^1 = (24,9,-14)$ is obtained and added to the representation $R(N)=\{z^1\}$. The new $\epsilon_3$ is updated according to the desired uniformity level as: $\epsilon_3 = z^{1}_3 + \delta_3 = -14+10 = -4$. Before the next iteration the worst value for objective function $z_2$, $z^{wv}_2$, is also updated to $9$.
    \item The second iteration, Figure \ref{fig:GPBA-B-2}, uses $\epsilon_3 =-4$ to solve Problem \ref{es-MOO} and point $z^2 = (24,5,-3)$ is obtained and added to the representation, $R(N)=\{z^1,z^2\}$. The third component of the $\epsilon$ vector is updated accordingly, becoming $\epsilon_3 = 7$. Since $z^2_2 < z^{wv}_2$, the worst value for objective function $z_2$ is updated, resulting in $z^{wv}_2 = 5$.
    \item When solving Problem \ref{es-MOO} with $\epsilon_2 = -28$ and $\epsilon_3 = 7$, Figure \ref{fig:GPBA-B-3}, point $z^3 = (18,8,9)$ is obtained and added to the representation, $R(N)=\{z^1,z^2,z^3\}$. The $\epsilon_3$ is updated as $\epsilon_3=z^{3}_3 + \delta_3 = 9+10 = 19$. Again, since $z^2_2\geqslant z^{wv}_2$, the worst value for objective function $z_2$ is not updated.
    \item In the fourth iteration, Figure \ref{fig:GPBA-B-4}, Problem \ref{es-MOO} is solved using $\epsilon_3=19$. The outcome vector is $z^4 = (12,11,21)$, which is added to the representation, $R(N)=\{z^1,z^2,z^3, z^4\}$. As a consequence, $\epsilon_3$ is updated as $z^{4}_3 + \delta_3 = 21+10 = 31$. In this iteration, because $z^2_2\geqslant z^{wv}_2$, the worst value for objective function $z_2$ does not need to be updated, remaining at $z^{wv}_2 = 5$.
    \item For the fifth iteration, Figure \ref{fig:GPBA-B-5}, Problem \ref{es-MOO} is solved using $\epsilon_3=31$. The resulting outcome vector is $z^5 = (6,14,33)$, which is added to the representation, becoming $R(N)=\{z^1,z^2,z^3, z^4, z^5\}$. Once more, $z^{wv}_2$ is not updated since $z^{wv}_2<14$. Then, $\epsilon_3$ is updated according to Algorithm \ref{alg:GPBA-B}, however, since the resulting $\epsilon_3$ is $43>z^{\ast}_3 = 42$ the loop is finished (see Figure \ref{fig:flow_generic}).
\end{enumerate}

For the next iteration over the outermost loop, objective function $z_2$, the value of $z^{wv}_2$ is used as $z_2$ to compute the next value of $\epsilon_2$, $\epsilon_2 = z_2 + \delta_2 = 5+10 = 15$.

\subsubsection{Cardinality representation (\textit{GPBA-C})}
\noindent Consider a desired cardinality for the representation of $25$ points, corresponding to a division of each constrained objective $k$, $k=1,\ldots,p,\; k\ne q$, into $5$. The first iteration over the outermost loop, $k=2$, starts with $\epsilon_2 = z^{nad}_2 = -28$ and $z^{wv}_2 = z^{\ast}_2 = 49$, and then iterates over the innermost loop in the following way:
\begin{enumerate}
    \item For the first iteration, in Figure \ref{fig:GPBA-C-1}, $\epsilon_3$ is $-48$, corresponding to the nadir approximation value for that objective, $z^{nad}_3$. Solving Problem \ref{es-MOO} with those parameters, yields to the outcome vector $z^1 = (24,9,-14)$, which is added to the representation $R(N)$, and we update $z^{wv}_2$ with $9$. Since this is the first iteration, the grid start point $z^{start}_3$ is $-48$, i.e.,  the nadir approximation value, the position within the grid, $i_3$, is $0$, and the number of solutions to obtain with such a grid, $c'_3$, is equal to the objective's cardinality minus one, $4$. The slack variable for Problem \ref{es-MOO}, $s_3$, is the difference between $\epsilon_3$ and the third component of the outcome vector, i.e., $s_3 = -14 - (-48) = 34$. To determine whether a step in the grid may be skipped or not, we firstly compute the step size: $step = \max\big\{(42-(-48))/4,\,1\big\} = 22.5$. Since the number of steps to skip, $b$, is the integer part of $34/22.5=1.51$, then $1$ step in the grid must be skipped. Consequently, the grid is refined from the obtained solution, $z^{start}_3 = z^1_3 = -14$, onwards, considering the number of points left to compute: $c'_3 = 4 - 0 = 4$, $step = \max\big\{(42-(-14))/4,\,1\big\} = 14$, $i_3 = 1$ (Figure \ref{fig:GPBA-C-2}). Hence, the $\epsilon_3$ parameter for the following iteration is $\epsilon_3 = -14 + 1\times 14 = 0$. 
    \item For the second iteration, in Figure \ref{fig:GPBA-C-2}, we use $\epsilon_3=0$ when solving Problem \ref{es-MOO}, leading to the vector $z^2 = (22,6,1)$, which is added to the representation, $R(N)=\{z^1,z^2\}$. The worst value for objective function $z_2$ needs to be updated with $6$, since $z_2^2<z^{wv}_2$. As a result, the 3\textsuperscript{rd} slack variable for Problem \ref{es-MOO} is $s_3 = 1 - 0=1$. Since the step size is $14$, no solution is skipped, and it can be determined by the integer part of $b = \lfloor 1/14\rfloor = 0$. Hence the grid remains the same and the position within the grid is increased by one, $i_3 = 2$. Hence, the $\epsilon_3$ parameter for the following iteration is $\epsilon_3 = -14 + 2\times 14 = 14$.
    \item For the third iteration, Figure \ref{fig:GPBA-C-3}, we solve Problem \ref{es-MOO} with $\epsilon_3 = 14$, leading to vector $z^3 = (14,13,14)$, which is added to the representation, $R(N)=\{z^1,z^2,z^3\}$. The worst value for objective function $z_2$, $z^{wv}_2$, remains the same. The slack variable for Problem \ref{es-MOO} is $s_3 = 0$, hence no solution are skipped and the grid is kept. The position within the grid is increased by one, $i_3 = 3$, and $\epsilon_3$ parameter is updated to $\epsilon_3 = -14 + 3\times 14 = 28$.
    \item For the forth iteration, Figure \ref{fig:GPBA-C-4}, we use $\epsilon_3 =28$, which results in the outcome vector $z^4 = (8,13,29)$, and the representation $R(N)=\{z^1,z^2,z^3,z^4\}$.  The worst value for objective function $z_2$ keeps its value $z^{wv}_2=6$, because $z^4_2\geqslant z^{wv}_2$. Since $z^4_3 - \epsilon_3 = 1$ corresponds to the slack variable, there are no redundant solutions in the grid, $b = \lfloor 1/14\rfloor = 0$. As a result the grid remains the same and the position within the grid is increased by one, $i_3 = 4$. The $\epsilon_3$ parameter is then updated to $\epsilon_3 = -14 + 4\times 14 = 42$.
    \item For the fifth and last iteration in this loop, Figure \ref{fig:GPBA-C-5}, we use $\epsilon_3 =42$, which results in the outcome vector $z^5 = (0,20,42)$, updating the representation to $R(N)=\{z^1,z^2,z^3,z^4,z^5\}$. Since $\epsilon_3 = z^5_3$ no solution was skipped and the position within the grid may be increased by one, $i_3 = 5$. Due to $i_3$ being greater than $c'_3$, the updated value of $\epsilon_3 = -14 + 5\times 14 = 56$ is greater than the ideal value $z^{ast}=42$. Hence, the loop will be exited and the grid control parameters for this objective, $k=3$, will return to the values of the first iteration.
\end{enumerate}

For the next iteration over the outermost loop, objective function $z_2$, the value of $z^{wv}_2=6$ is used as $z_2$ to compute the next value of $\epsilon_2$, using Algorithm \ref{alg:GPBA-C}. Because the slack variable $s_2 = -28-6=-34$ and the step size to $step = \max\big\{(49-(-28))/4,\,1\big\} = 19.25$, resulting in skipping one step, $b=\lfloor|-34/19.25|\rfloor=1$. The grid parameters are adapted to $z^{start}_2 = z_2 = 6$, $c'_2 = 4-0 = 4$ and $i_2 = 1$. The step size is updated to $step = \max\big\{(49-6))/4,\,1\big\} = 10.75$ and, at last, $\epsilon_2 = 6 + 1\times 10.75 = 16.75$.

\subsection{Generating the entire Pareto front}
\noindent The three algorithms presented in the above section, \textit{GPBA-A}, \textit{GPBA-B}, and \textit{GPBA-C}, are based on solving a sequence of \ref{es-MOO} problems, hence they can be used to generate the complete Pareto front under the condition that the objective function only takes integer values. Additionally the parameters that control each algorithm, the acceptable coverage error, the uniformity level, and the cardinality, must be chosen to ensure the production of the entire Pareto set. 

When using algorithm \textit{GPBA-A} to compute the entire Pareto front, the desired acceptable coverage error, $\gamma$, is $1$, guaranteeing that the distance between any two consecutive points in the Pareto front representation is at most unitary. Since this algorithms are applied to problems of type \ref{MOO}, if $\Gamma \leqslant 1$, then $\Gamma$ is, in fact, equal to $1$.

In case algorithm \textit{GPBA-B} is used, the acceptable uniform level, $\delta$, should be set to $1$ to achieve the computation of the entire Pareto front. Setting $\delta=1$ ensures the algorithm will search for the next solution as close to the previous as a unitary step in the objective space. 

At last, for the case of \textit{GPBA-C}, when generating the entire Pareto front, the cardinality should be set for each objective as the respective range. As a result, a unitary step size is used for each objective function's loop and the algorithm restricts the search region for the next solution to a unitary distance in the objective space.

All three proposed algorithms will be tested in the computational experiments section generating the entire Pareto front, hence using the parameters setting described above. 


\section{Computational experiments}
\noindent This section is devoted to present the design of the experiments, some implementation issues, the computational results for the identification of the whole Pareto front, some experiments on the representations of the Pareto front, and a few final comments and remarks.

The performance of the three proposed algorithms \textit{GPBA-A}, \textit{GPBA-B} and \textit{GPBA-C} is compared with the works of \cite{MavrotasFlorios2013} and \cite{ZhangReimann2014}, hereinafter referred to as \textit{AUGM-2} and \textit{S-AUGM}, respectively. All algorithms were coded using PyCharm Community Edition 2020 for Windows Desktop, Python 3.7 as interpreter and CPLEX 12.9 as optimization solver. All experiments were performed in a workstation equipped with two processors Intel Xeon X5680 3.33GHz and 24GB of RAM, running Windows 10.

\subsection{The design of the experiments}
\noindent  The generation of the instances used in these experiments was performed through a generator developed by the second author of this study\footnote{For more details please contact José Rui Figueira at figueira@tecnico.ulisboa.pt}. This generator makes use of the random number generation procedure of NETGEN generator proposed in \cite{KlingmanEtAl1974} and requires as input the following elements:
\begin{itemize}[label={--}]
    \item The number of variables (\textit{n}), objective functions (\textit{p}), and constraints (\textit{m}).
    \item The ranges for the parameter values (coefficients of the constraints and the weights/profits of the objective functions). These parameters are randomly generated by using a uniform distribution with a specific seed number for each constraint/objective function. 
    \item In this study, all seed numbers are increased by a fixed value, for each instance generated. 
    \item The right-hand-side of each constraint is bounded by half the sum of the constraint's coefficients.
\end{itemize}

The algorithm's performance is essentially compared in terms of solution quality and running time, or a proxy of the latter, such as the number of iterations \textit{per} supported outcome vector computed. As such, since both running time and its proxies are often not normally distributed, each problem type was tested on 30 instances. The choice of 30 instances was motivated by the conclusion of \cite{CoffinSaltzman2000} that, for those measures of performance, 30 to 40 observation would guarantee statistically significant results. Apart from the aforementioned observation, \cite{CoffinSaltzman2000} also provided other suggestions to draw meaningful conclusions from the comparison of algorithms that were followed in this work, such as the use of sets of problems where certain parameters are varied, and the use of mean and variance of metrics for comparison. As a result, a total of 570 different uncorrelated instances were generated.

\subsection{Implementation issues}\label{sec: implementation issues}
\noindent All algorithms, aside from \textit{AUGM-2}, can overcome poor nadir estimations, since the constraints' right-hand-side is computed based on the obtained solution and not beforehand. Hence, for those algorithms, in these experiments, we simply use the single minimization of each objective function individually as a nadir point. However, since the improvements on the constrained objective functions are weighted by the amplitude of the objective's range, due to some solvers' sensitivity, the value of parameter $\rho$ in Problem \ref{es-MOO} may have to be adapted for some instances. The consequence of not adapting parameter $\rho$ is the failure to compute some non-dominated solutions.

Algorithm \textit{AUGM-2}, however, requires good quality nadir estimations. For such a purpose, as suggested in \cite{MavrotasFlorios2013}, the lexicographic pay-off table is used to provide an approximation of the nadir point. However, as it has been discussed in literature, lexicographic pay-off tables do not guarantee a good estimation of nadir points for more than two objective functions \citep[see, e.g.,][]{IsermannSteuer1988,AlvesClimaco2007}. For this implementation reason, in some instances, \textit{AUGM-2} does not explore the full Pareto front and, consequently, skips non-dominated solutions.

\subsection{Experiments on generating the whole Pareto front}
\noindent To evaluate the computational and accuracy behaviour of different algorithms, when computing the complete Pareto front, experiments will be performed on binary and integer instances, considering both single and multi-dimensional cases. Additionally, $180$ tests are run on multi-objective general integer programming instances. This is done in order to, not only compare the algorithms, but understand how changes in the number of objective functions, constraints, variables, and problem types affect the performance of the algorithms. 

In the remaining of this section, analysis of the results for the aforementioned experiments are presented. This analysis is accompanied by summary tables, where underlined are results that show algorithms' drawbacks, while in bold are results that show algorithms' strengths. Algorithms are evaluated on the average and standard deviation over the number of unique non-dominated solutions found ($\vert \overline{N(Z)} \vert$ and $\sigma_{N(Z)}$) and  CPU time, as well as the average number of iterations computed for each non-dominated solution found ($\overline{iter/\hat{z}}$). Although the latter corresponds to a proxy of the CPU time, it is included in this study since it is not correlated to the number of solutions in the Pareto front while the CPU time typically is, better conveying the influence of the remaining problem parameters on the algorithm's performance.

\subsubsection{Multi-objective $\{0,1\}-$knapsack instances}\label{sec: Multi-objective 0,1-knapsack instances}
\noindent The first set of problems studied were multi-objective binary knapsack instances with three and four objective functions, considering one constraint. The three-objective instances have $50$, $75$, and $100$ variables, while the four-objective instances were generated with just $50$ variables. Both the objective function's profits and the constraint's coefficients are drawn from a uniform distribution of the type $U[1,100]$. The initial seed numbers for the three-objective instances with $50$ variables were $128$, $888$ and $6$ for each objective function profit, and $40$ for the constraint's coefficients. For the remaining three-objective instances, the seed numbers for the objective functions' profits were $47$, $28$ and $626$, respectively, whilst the seed number $135$ was used for the constraint's coefficients. At last, for the four-objective instances, the seed numbers of the objective functions' profits are $47$, $28$, $626$ and $135$, and $298$ for the constraint's coefficients. All seed numbers are increased by $5$ on each instance generated of each type. For each set of problems, $30$ instances were generated. Table \ref{table: 0-1 knapsack - 1 const} presents experiments results for the above described instances.

\begin{table}[!htbp]
\centering
\begin{tabular}{rrlrrrrr} 
 \hline\hline
 $p$ & $n$ & Algorithm & $\vert \overline{N(Z)} \vert$ & $\sigma_{N(Z)}$ & $\overline{iter/\hat{z}}$  & $CPU$ (sec) & $\sigma_{CPU}$ \TBstrut\\[0.5ex] 
 \hline\hline
 3 & 50 & \textit{AUGM-2}  &  \underline{403.47}   &  209.13   &  \underline{78.76} &  \underline{1133.90}  &  773.31   \\
   &    & \textit{S-AUGM}  &  408.27   &  210.30   &  $\mathbf{1.89}$  &   58.36    &  34.51    \\
   &    & \textit{GPBA-A}  &  408.27   &  210.30   &  \underline{14.37} &  \underline{394.56}   &  238.80   \\
   &    & \textit{GPBA-B}  &  408.27   &  210.30   &   $\mathbf{1.89}$  &  $\mathbf{58.36}$    &  $\mathbf{34.04}$    \\
   &    & \textit{GPBA-C}  &  408.27   &  210.30   &   $\mathbf{1.89}$  &  58.61   &  34.64    \\\hline
   & 75 & \textit{AUGM-2}  &  \underline{1512.63}  &  915.60   &  \underline{52.09} &  \underline{3628.31}  &  1885.52  \\
   &    & \textit{S-AUGM}  &  1518.83  &  920.27   &   $\mathbf{1.82}$  &  263.34   &  193.64   \\
   &    & \textit{GPBA-A}  &  1518.83  &  920.27   & \underline{11.69}  &  \underline{1771.25}  &  1378.51  \\
   &    & \textit{GPBA-B}  &  1518.83  &  920.27   &   $1.83$  &   $\mathbf{260.03}$   &  190.42   \\
   &    & \textit{GPBA-C}  &  1518.83  &  920.27   &   $1.83$  &  260.18   &  189.67   \\\hline 
   &100 & \textit{AUGM-2}  &  \underline{3343.37}  &  1609.82  & \underline{42.82}  &  \underline{9237.81} &  4656.32 \\
   &    & \textit{S-AUGM}  &  3356.40  &  1628.35  & \textbf{1.79}   &  666.86 &  411.67 \\
   &    & \textit{GPBA-A}  &  3356.40  &  1628.35  & \underline{9.28}  & \underline{4501.94} &  2674.97  \\
   &    & \textit{GPBA-B}  &  3356.40  &  1628.35  & \textbf{1.79}   & \textbf{661.26} &   427.21  \\
   &    & \textit{GPBA-C}  &  3356.40  &  1628.35  & \textbf{1.79}   &  662.42 &  427.22 \\ \hline
  4  & 50 & \textit{S-AUGM}	&	2655.43	&	1392.73	&	$\mathbf{4.74}$	&	\textbf{9940.89}		&	7008.43 \\
     &    & \textit{GPBA-B} &	2655.43	&	1392.73	&	$\mathbf{4.74}$	&	\underline{10104.69}	&	7152.10\\
     &    & \textit{GPBA-C} & 2655.43	&	1392.73	&	\underline{4.77}	&	10032.92   & 7205.45\\
 \hline\hline
\end{tabular}
\caption{Multi-objective 0-1 knapsack instances}
\label{table: 0-1 knapsack - 1 const}
\end{table}

Regarding the instances with three objectives, the computational effort, in terms of computational time, increases significantly as the number of solutions in the Pareto front and problem variables also increases. However, the increase in the number of variables does not impact on the number of iterations \textit{per} non-dominated solution computed. All algorithms apart from \textit{AUGM-2}, were able to compute the full Pareto front. Algorithm \textit{AUGM-2} failed to obtain some non-dominated solutions due to the strategy used to provide the nadir point, which can generate poor quality nadir estimations (refer to Section \ref{sec: implementation issues}, for more details). Additionally, \textit{AUGM-2} is the least efficient method among the ones studied, both in terms of computational time and in terms of number of iterations \textit{per} non-dominated solution obtained. The reason for this is that, in this algorithm, iteration skipping and early exit mechanisms are only applied to the innermost loop. \textit{GPBA-A} also does not present competitive results since more iterations are performed \textit{per} non-dominated solution computed for this algorithm when compared to the others, which is also reflected in the computational time. \textit{GPBA-A} presents this behaviour due to the fact that the accelerated early exit strategy cannot be applied to this algorithm. The results displayed by algorithms \textit{AUGM-2} and \textit{GPBA-A} justify the need and added value of the integration of the acceleration strategies (Figures \ref{fig:flow_check_redundancy} and \ref{fig:early_exit}). Consequently, algorithms \textit{AUGM-2} and \textit{GPBA-A} were discarded from further analysis. At last, when comparing the performance of \textit{S-AUGM} and the other two proposed algorithms, \textit{GPBA-B} and \textit{GPBA-C}, the presented results are similar and all appear to be competitive.  

In the case of the four-objective instances, for the aforementioned reasons, only algorithms \textit{GPBA-B} and \textit{GPBA-C} are compared with \textit{S-AUGM}. From the analysis of the results it is possible to observe that the three compared algorithms show a high sensitivity to the increase in the number of objectives, having increased the average number of iterations needed to compute a non-dominated solution. Furthermore, although all algorithms remained competitive, \textit{GPBA-C} appears to have the lowest efficiency having performed, on average, the most iterations per non-dominated solution found. This could be explained by possibly needing a readjustment of parameter $\rho$ in Problem \ref{es-MOO}, as explained in Section \ref{sec: implementation issues}. However, when comparing in terms of computational time, \textit{GPBA-C} is slightly more efficient than \textit{GPBA-B}.

\subsubsection{Multi-objective multi-dimensional $\{0,1\}-$knapsack instances}\label{sec: Multi-objective multi-dimensional 0-1}
\noindent For the multi-objective multi-dimensional binary knapsack problems, instances were generated varying the number of constraints between $2$ and $5$, with a fixed number of objectives at $3$ and variables at $25$. Both the objective function's profits and the constraint's coefficients follow a uniform distribution bounded between 1 and 100. The seed numbers used for the objective function's profits were $47$, $63$ and $728$, while for each constraint equation's coefficient the values used were 626, 135, 28, 17, 758, respectively. Table \ref{table: 0-1 knapsack - multi-dimendion} displays the results for algorithms \textit{S-AUGM}, \textit{GPBA-B} and \textit{GPBA-C}. For the reasons conveyed in Section \ref{sec: Multi-objective 0,1-knapsack instances}, related to poor performance, algorithms \textit{AUGM-2} and \textit{GPBA-A} are not considered in this comparison.


\begin{table}[!htbp]
\centering
\begin{tabular}{rrrlrrrrr} 
 \hline\hline
 $p$ & $n$ &  $m$ &Algorithm & $\vert \overline{N(Z)} \vert$ & $\sigma_{N(Z)}$ & $\overline{iter/\hat{z}}$  & $CPU$ (sec) & $\sigma_{CPU}$ \TBstrut\\[0.5ex]
 \hline\hline
  3  & 25 & 2  & \textit{S-AUGM} &	80.03	&	38.39	&	$\mathbf{1.91}$	&	$\mathbf{8.06}$	&	4.53 \\
     &    &    & \textit{GBPA-B} &	80.03	&	38.39	&	$\mathbf{1.91}$	&	8.11	&	4.49 \\
     &    &    & \textit{GBPA-C} &	80.03	&	38.39	&	$\mathbf{1.91}$	&	8.30	&	4.74 \\\hline
     &    & 3  & \textit{S-AUGM} &	84.83	&	40.56	&	$\mathbf{1.92}$	&	$\mathbf{8.95}$	&	4.81 \\
     &    &    & \textit{GPBA-B} &	84.83	&	40.56	&	$\mathbf{1.92}$	&	8.99	&	4.95 \\
     &    &    & \textit{GPBA-C} &	84.83	&	40.56	&	$\mathbf{1.92}$	&	9.03	&	4.90 \\ \hline 
     &    & 4  & \textit{S-AUGM} &	88.30	&	44.96	&	$\mathbf{1.92}$	&	$\mathbf{9.67}$	&	5.49 \\
     &    &    & \textit{GPBA-B} &	88.30	&	44.96	&	$\mathbf{1.92}$	&	9.78	&	5.58 \\
     &    &    & \textit{GPBA-C} &	88.30	&	44.96	&	$\mathbf{1.92}$	&	9.73	&	5.56 \\ \hline
     &    & 5  & \textit{S-AUGM} &	90.73	&	37.97	&	$\mathbf{1.92}$	&	$\mathbf{9.94}$	&	4.71 \\
     &    &    & \textit{GPBA-B} &	90.73	&	37.97	&	$\mathbf{1.92}$	&	10.02	&	4.85 \\
     &    &    & \textit{GPBA-C} &	90.73	&	37.97	&	$\mathbf{1.92}$	&	10.12	&	4.89 \\
 \hline\hline
\end{tabular}
\caption{Multi-objective multi-dimensional 0-1 knapsack instances}
\label{table: 0-1 knapsack - multi-dimendion}
\end{table}

It is possible to observe that the performance of the algorithms is not affected by the increase in the number of constraints. Additionally, the three tested algorithms, \textit{S-AUGM}, \textit{GPBA-B} and \textit{GPBA-C}, are competitive, finding the full Pareto front in a short amount of time. \textit{S-AUGM} appears to be slightly more efficient when comparing average computational time, however differences between algorithms represent, in the worst case, $1.8\%$.

\subsubsection{Multi-objective integer knapsack instances}
\noindent To test the algorithms on integer multi-objective knapsack instances, the instances generated in Section \ref{sec: Multi-objective 0,1-knapsack instances} with three objectives and 50 variables were solved using integer variables. Instance number 27 that uses seed numbers 158, 1018 and 136 for the objective function's profits and 170 for the constraint's coefficients, was discarded since, by having many Pareto front solutions, was taking too long to solve. Hence, the results presented in Table \ref{table: Multi-objective integer knapsack instances} consider 29 instances. 


\begin{table}[!htbp]
\centering
\begin{tabular}{rrlrrrrr} 
 \hline\hline
 $p$ & $n$ & Algorithm & $\vert \overline{N(Z)} \vert$ & $\sigma_{N(Z)}$ & $\overline{iter/\hat{z}}$  & $CPU$ (sec) & $\sigma_{CPU}$ \TBstrut\\[0.5ex]
 \hline\hline
 3	&	50	&	\textit{S-AUGM}	&	3152.59	&	9061.66	&	\underline{1.98}    &	$\mathbf{2590.92}$	&	6666.79	\\
	&		&	\textit{GPBA-B}	&	3152.59	&	9061.66	&	$\mathbf{1.71}$	&	$2662.00$	&	7017.92	\\
	&		&	\textit{GPBA-C}	&	3152.59	&	9061.66	&	$\mathbf{1.71}$	&	$2673.95$	&	7051.64	\\
 \hline\hline
\end{tabular}
\caption{Multi-objective integer knapsack instances}
\label{table: Multi-objective integer knapsack instances}
\end{table}

Although all algorithms are able to compute the full Pareto front, the performance differs. If on one hand \textit{S-AUGM} takes less computational time than the other two algorithms, around $3\%$ less than \textit{GPBA-C}, it is also the least efficient method when assessing the number of iteration needed to compute each non-dominated solution. This suggests that, both \textit{GPBA-B} and \textit{GPBA-C} computational time performance could be bounded by code efficiency.

\subsubsection{Multi-objective multi-dimensional integer knapsack instances}\label{sec: Multi-objective multi-dimensional integer knapsack instances}
\noindent To assess the impact of the increase in dimensions on the algorithms performance on integer instances, the instances from Section \ref{sec: Multi-objective multi-dimensional 0-1} were applied using integer type variables. Table \ref{table: Integer multi-dimensional} presents the results for this case.

\begin{table}[!htbp]
\centering
\begin{tabular}{rrrlrrrrr} 
 \hline\hline
 $p$ & $n$ &  $m$ &Algorithm & $\vert \overline{N(Z)} \vert$ & $\sigma_{N(Z)}$ & $\overline{iter/\hat{z}}$  & $CPU$ (sec) & $\sigma_{CPU}$ \TBstrut\\[0.5ex]
 \hline\hline
  3	&	25	&	2	&	\textit{S-AUGM} &	167.07	&	269.14	&	\underline{1.93}	&	\underline{30.13}	&	66.06\\
	&		&		&	\textit{GPBA-B} &	167.07	&	269.14	&	$\mathbf{1.89}$	&	$\mathbf{28.72}$	&	63.53\\
	&		&		&	\textit{GPBA-C} &	167.07	&	269.14	&	$\mathbf{1.89}$	&	28.79	&	62.97\\ \hline
	&		&	3	&	\textit{S-AUGM} &	292.57	&	302.10	&	\underline{1.94}	&	\underline{50.94}	&	58.59\\
	&		&		&	\textit{GPBA-B} &	292.57	&	302.10	&	$\mathbf{1.88}$	& $\mathbf{28.72}$	&	54.10\\
	&		&		&	\textit{GPBA-C} &	292.57	&	302.10	&	$\mathbf{1.88}$	&	28.79	&	53.83\\ \hline
	&		&	4	&	\textit{S-AUGM} &	285.50	&	232.03	&	\underline{1.93}	&	\underline{53.51}	&	44.51\\
	&		&		&	\textit{GPBA-B} &	285.50	&	232.03	&	$\mathbf{1.92}$	&	$\mathbf{49.55}$	&	44.51\\
	&		&		&	\textit{GPBA-C} &	285.50	&	232.03	&	$\mathbf{1.92}$	&	49.66	&	44.33\\ \hline
	&		&	5	&	\textit{S-AUGM} &	238.30	&	169.77	&	$\mathbf{1.94}$	&	\underline{42.30}	&	32.86\\
	&		&		&	\textit{GPBA-B} &	238.30	&	169.77	&	$\mathbf{1.94}$	&	$\mathbf{40.04}$	&	31.53\\
	&		&		&	\textit{GPBA-C} &	238.30	&	169.77	&	$\mathbf{1.94}$	&	40.08	&	31.96\\ 
 \hline\hline
\end{tabular}
\caption{Multi-objective multi-dimensional integer knapsack instances}
\label{table: Integer multi-dimensional}
\end{table}

Although the impact of the increase in dimension does not appear to be significant, the same trend as in Table \ref{table: Multi-objective integer knapsack instances} can be observed where \textit{S-AUGM} is less efficient than \textit{GPBA-B} and \textit{GPBA-C}. In this case, \textit{GPBA-B} and \textit{GPBA-C} outperform consistently \textit{S-AUGM}, both in terms of number of iterations solved \textit{per} non-dominated solution obtained and in terms of computational time. The latter indicator also suggests that \textit{GPBA-B} is more efficient than \textit{GPBA-C}. The worse performance of \textit{S-AUGM} in integer instances, when compared to the other algorithms, can probably be attributed to the $\epsilon$-constraint model formulation used on \textit{S-AUGM} not using elastic constraints as proposed by \cite{EhrgottRyan2003}. However, it is possible to observe on Table \ref{table: Integer multi-dimensional} that, as the number of constraints increase, the performance difference between \textit{S-AUGM} and the other algorithms appears to diminish.

\subsubsection{Multi-objective general integer programming instances}
\noindent Instances for general integer problems were generated, setting a uniform distribution ranging between 0 and 20 for the objective function's profits and the constraints' coefficients. The instances were generated for both three and four objectives with multiple constraints and 10, 15 and 20 variables. The generation of objective function's profits started with seed numbers 47, 63, 728 and 11, one for each objective function. The initial seed numbers used to generate the constraints' coefficients was 634, 17, 28, 626, 135, 34, 78, 55, 783, 945, 823, 362, 133, 92 and 41, one for each constraint equation. Seed numbers were increased by 5 on each instance generated of each type. Results on solving these instances using \textit{S-AUGM}, \textit{GPBA-B} and \textit{GPBA-C} are displayed in Table \ref{table: GenInteger}

\begin{table}[!htbp]
\centering
\begin{tabular}{rrrlrrrrr} 
 \hline\hline
 $p$ & $n$ &  $m$ & Algorithm & $\vert \overline{N(Z)} \vert$ & $\sigma_{N(Z)}$ & $\overline{iter/\hat{z}}$  & $CPU$ (sec) & $\sigma_{CPU}$ \TBstrut\\[0.5ex]
 \hline\hline
  3    &    10   &    5    &    \textit{S-AUGM}     &    19.83     &    10.14     &    $\mathbf{1.80}$     &    \underline{1.54}    &    1.06\\
     &         &         &    \textit{GPBA-B}     &    19.83     &    10.14     &    $\mathbf{1.80}$     &    $\mathbf{1.37}$     &    0.64\\
     &         &         &    \textit{GPBA-C}     &    19.83     &    10.14     &    $\mathbf{1.80}$     &    1.45                &    0.75\\\hline
     &    15   &    10   &    \textit{S-AUGM}     &    38.93     &    20.13     &    \underline{1.79}    &    \underline{3.72}    &    2.16 \\
     &         &         &    \textit{GPBA-B}     &    38.93     &    20.13     &    $\mathbf{1.78}$     &    $\mathbf{3.36}$     &    1.97\\
     &         &         &    \textit{GPBA-C}     &    38.93     &    20.13     &    $\mathbf{1.78}$     &    3.45                &    1.96\\\hline
     &    20   &    15   &    \textit{S-AUGM}     &    75.13     &    31.07     &    $\mathbf{1.73}$     &    \underline{9.03}    &    3.91 \\
     &         &         &    \textit{GPBA-B}     &    75.13     &    31.07     &    $\mathbf{1.73}$     &    $\mathbf{8.74}$     &    3.95\\
     &         &         &    \textit{GPBA-C}     &    75.13     &    31.07     &    $\mathbf{1.73}$     &    8.91                &    3.85\\\hline
4    &    10   &    5    &    \textit{S-AUGM}     &    30.67     &    13.31     &    $\mathbf{3.54}$     &    $\mathbf{7.73}$     &    4.71\\
     &         &         &    \textit{GPBA-B}     &    30.67     &    13.31     &    \underline{3.58}    &    7.75                &    4.76\\
     &         &         &    \textit{GPBA-C}     &    30.67     &    13.31     &    3.55                &    \underline{7.86}    &    4.75\\\hline
     &    15   &    10   &    \textit{S-AUGM}     &    107.77    &    74.05     &    $\mathbf{3.60}$     &    $\mathbf{43.39}$    &    37.44\\
     &         &         &    \textit{GPBA-B}     &    107.77    &    74.05     &    \underline{3.64}    &    44.23               &    39.23\\
     &         &         &    \textit{GPBA-C}     &    107.77    &    74.05     &    3.63                &    \underline{44.70}   &    38.86\\\hline
     &    20   &    15   &    \textit{S-AUGM}     &    311.30    &    175.06    &    $\mathbf{3.52}$     &    $\mathbf{152.01}$   &    94.09\\
     &         &         &    \textit{GPBA-B}     &    311.30    &    175.06    &    \underline{3.56}    &    156.39              &    97.36\\
     &         &         &    \textit{GPBA-C}     &    311.30    &    175.06    &    3.55                &    \underline{156.43}  &    97.56\\ 
 \hline\hline
\end{tabular}
\caption{General multi-objective integer linear instances}
\label{table: GenInteger}
\end{table}

Regarding the results on the three objective instances, it is possible to observe the same trend as in Section \ref{sec: Multi-objective multi-dimensional integer knapsack instances} in which \textit{GPBA-B} and \textit{GPBA-C} outperform \textit{S-AUGM}, although not as significantly. The less significant difference in performance between \textit{S-AUGM} and the remaining algorithms may be attributed to these problems having more constraints than the ones presented in the previous section, where it was observed that the performance gap was lower as the number of constraints increased. Finally, \textit{GPBA-B} is consistently the best performing algorithm among the three.

When looking at the four objective instances, the same trend as before is not observed. On the contrary, \textit{S-AUGM} performs better than \textit{GPBA-B} and \textit{GPBA-C} in all sets of problems, both in terms of computational time and in terms of number of iterations \textit{per} non-dominated solution computed. \textit{GPBA-B}, for all sets of problems, requires the most iterations \textit{per} non-dominated solution, while \textit{GPBA-C} takes the most time. 

\subsection{Experiments on representation}
\noindent To assess the quality of the proposed algorithms for the representation problem, experiments were carried out on the binary knapsack problem instances, Section \ref{sec: Multi-objective 0,1-knapsack instances}, varying the target cardinality for each objective ($c_k$, for all $k = {1,\dots,p}$, $k\neq q$). As a pre-computing stage for each instance, the acceptable uniformity level, $\delta_k$, and the acceptable coverage level, $\gamma_k$, are defined as $(z^*_k - z^{nad}_k)/c_k$, where the estimation of the nadir point $z^{nad}_k$ was, only for this purpose of defining the acceptable uniformity and coverage level, obtained using the lexicographic pay-off table. Table \ref{table: Representation} presents the results for these experiments. Apart from the already introduced performance measures, the obtained representations are evaluated in terms of the  coverage error, $\Gamma\big(R(N)\big)$, the uniformity level, $\Delta\big(R(N)\big)$, and the cardinality, $\Pi\big(R(N)\big)$. In this section, only the three proposed algorithms are compared since \textit{S-AUGM} was developed exclusively to compute the whole Pareto front and not a representation of it.

\begin{table}[!htbp]
\centering
\begin{tabular}{rrrlrrrrr} 
\hline\hline
$p$ & $n$ & $c_k$ & Algorithm & $\vert\Pi\big(R(N)\big)\vert$ & $\vert\Gamma\big(R(N)\big)\vert$ & $\vert\Delta\big(R(N)\big)\vert$ & $\overline{iter/\hat{z}}$ & $CPU$(s)\TBstrut\\[0.5ex] \hline\hline
3 & 50 & 5 & \textit{GPBA-A} & \underline{25.83} & \textbf{186.67} & \underline{29.43} & \underline{1.42} & \underline{1.45} \\
 &  &  & \textit{GPBA-B} & \textbf{13.60} & \underline{274.33} & \textbf{46.53} & \textbf{1.16} & \textbf{0.78} \\
 &  &  & \textit{GPBA-C} & 15.33 & 205.00 & 42.93 & 1.22 & 0.89 \\ \hline
 &  & 17 & \textit{GPBA-A} & \underline{113.70} & \textbf{113.53} & \underline{8.20} & \underline{1.80} & \underline{8.64} \\
 &  &  & \textit{GPBA-B} & \textbf{67.40} & \underline{160.90} & \textbf{14.40} & \textbf{1.20} & \textbf{3.65} \\
 &  &  & \textit{GPBA-C} & 112.20 & 114.53 & 8.30 & 1.36 & 6.86 \\ \hline
 &  & 25 & \textit{GPBA-A} & 153.53 & 98.77 & 6.70 & \underline{2.13} & \underline{13.80} \\
 &  &  & \textit{GPBA-B} & \textbf{102.17} & \underline{143.63} & \textbf{10.03} & \textbf{1.25} & \textbf{5.71} \\
 &  &  & \textit{GPBA-C} & \underline{168.83} & \textbf{94.13} & \underline{5.87} & 1.57 & 11.56 \\ \hline
3 & 75 & 5 & \textit{GPBA-A} & \underline{31.50} & \textbf{279.17} & \underline{27.03} & \underline{1.33} & 1.99 \\
 &  &  & \textit{GPBA-B} & \textbf{15.20} & \underline{425.13} & 50.27 & \textbf{1.15} & \textbf{0.97} \\
 &  &  & \textit{GPBA-C} & 16.07 & 313.57 & \textbf{52.07} & 1.19 & \textbf{0.97} \\ \hline
 &  & 17 & \textit{GPBA-A} & \underline{199.33} & \textbf{157.10} & \underline{6.80} & \underline{1.41} & \underline{14.35} \\
 &  &  & \textit{GPBA-B} & \textbf{99.90} & \underline{214.50} & \textbf{12.30} & \textbf{1.14} & \textbf{6.05} \\
 &  &  & \textit{GPBA-C} & 152.13 & 175.33 & 8.40 & 1.18 & 9.50 \\ \hline
 &  & 25 & \textit{GPBA-A} & \underline{284.73} & \textbf{137.10} & \underline{4.97} & \underline{1.57} & \underline{22.65} \\
 &  &  & \textit{GPBA-B} & \textbf{169.17} & \underline{186.13} & \textbf{8.27} & \textbf{1.16} & \textbf{10.26} \\
 &  &  & \textit{GPBA-C} & 274.93 & 139.73 & 5.23 & 1.26 & 18.31 \\ \hline
3 & 100 & 5 & \textit{GPBA-A} & \underline{32.17} & \textbf{334.40} & \underline{29.00} & \underline{1.31} & 2.25 \\
 &  &  & \textit{GPBA-B} & \textbf{15.13} & \underline{530.47} & \textbf{87.00} & \textbf{1.17} & \textbf{1.09} \\
 &  &  & \textit{GPBA-C} & 16.17 & 397.20 & 55.07 & 1.20 & \textbf{1.09} \\ \hline
 &  & 17 & \textit{GPBA-A} & \underline{232.10} & \textbf{189.33} & \underline{6.80} & \underline{1.27} & \underline{17.11} \\
 &  &  & \textit{GPBA-B} & \textbf{109.90} & \underline{277.87} & \textbf{11.47} & \textbf{1.12} & \textbf{7.59} \\
 &  &  & \textit{GPBA-C} & 167.73 & 217.97 & 8.03 & \textbf{1.12} & 11.53 \\ \hline
 &  & 25 & \textit{GPBA-A} & \underline{361.07} & \textbf{161.70} & \underline{5.03} & \underline{1.37} & \underline{28.90} \\
 &  &  & \textit{GPBA-B} & \textbf{198.23} & \underline{227.50} & \textbf{10.20} & \textbf{1.12} & \textbf{13.43} \\
 &  &  & \textit{GPBA-C} & 323.97 & 174.73 & 5.27 & 1.16 & 23.45 \\
\hline\hline
\end{tabular}
\caption{Representation for multi-objective 0-1 knapsack instances}
\label{table: Representation}
\end{table}

On the one hand, in Table \ref{table: Representation} it is possible to observe that, apart from few exceptions, \textit{GPBA-A} presents the best results on coverage error and \textit{GPBA-B} on uniformity level. However, both algorithms show the worst performance in the other measure, i.e. \textit{GPBA-A} on uniformity level and \textit{GPBA-B} on coverage error. On the other hand, algorithm \textit{GPBA-C} proves to be an algorithm with a more balanced performance. When low values for the target cardinality are set, \textit{GPBA-C} obtains a good uniformity level, while not compromising coverage as much as \textit{GPBA-B}. For higher values of target cardinality, \textit{GPBA-C} behaves more like \textit{GPBA-A}, presenting low coverage errors and worsening the results on the uniformity level. Nevertheless, \textit{GPBA-C} is much more efficient than \textit{GPBA-A}, as can be observed by the number of models solved \textit{per} non-dominated solution in the representation, and the computational time. 

\subsection{Final comments and remarks}
\noindent All the three algorithms were tested on experiments set for computing the whole Pareto front and a representation of it. Regarding the computation of the whole Pareto front:
\begin{enumerate}
    \item All the three algorithms outperformed \textit{AUGM-2} algorithm. The reason for this is that \textit{AUGM-2} only skips redundant solution and applies the early exit mechanism in the innermost loop. Furthermore, the redundant solution skip strategy is less sophisticated than the one applied on the other algorithms since it does not consider all past solutions, but just the last one. These factors imply that \textit{AUGM-2} requires more models to be solved \textit{per} Pareto front solution computed.
    \item Algorithms \textit{GPBA-B} and \textit{GPBA-C} outperform \textit{GPBA-A}. \textit{GPBA-A}, due to the strategy used to determine the $\epsilon$ vector, does not have the early exit mechanism, lagging behind the other algorithms in terms of number of models solved \textit{per} Pareto front solution computed.
    \item Algorithms \textit{GPBA-B} and \textit{GPBA-C} are competitive with respect to \textit{S-AUGM}, being amongst the best in the literature. \textit{GPBA-B} and \textit{GPBA-C} showed very good results, similar to the ones presented by \textit{S-AUGM}. This was to be expected since the acceleration mechanisms are common to the three algorithms. However, probably due to the model structure used, \textit{S-AUGM} falls behind on integer knapsack instances with a low number of constraints, less than four. But, for general multi-objective integer instances, \textit{S-AUGM} appears to be more efficient.
\end{enumerate}

\noindent Regarding the representation of the Pareto front, the three algorithms were not compared with others in literature since, to the best of our knowledge, this are the first $\epsilon$-constraint based algorithms to target the representation problem for integer problems with more than two objectives. Furthermore, \textit{S-AUGM} did not present any extension for generating a representation of the Pareto front and \textit{AUGM-2} was inefficient in its search strategy. Hence, the analysis of the proposed algorithms performance on the representation of Pareto front showed that:
\begin{enumerate}
    \item \textit{GPBA-A} was the algorithm that presented the best coverage error, while having the worst uniformity level.
    \item \textit{GPBA-B} had the best uniformity level, presenting the worst coverage error.
    \item Contrary to what is observed for \textit{GPBA-A} and \textit{GPBA-B}, \textit{GPBA-C} behaviour changes depending on the target cardinality. When low cardinality is used, \textit{GPBA-C} appears to privilege uniformity behaving as \textit{GPBA-B}. When higher cardinality values are in place, \textit{GPBA-C} privileges coverage, behaving as \textit{GPBA-A} although more efficiently.
\end{enumerate}
\section{Conclusions}
\noindent This work addresses the Pareto front representation problem, contributing with methodologies to explore the feasible region of the objective space. This is relevant for most multi-objective problems, since the evaluation of the whole Pareto front may be overwhelming for the decision-maker, possibly inhibiting a decision. A good representation of the Pareto front is one that consists of as few solutions as possible, covers all regions of the objective space and spreads the solutions the most.  To tackle this problem three algorithms are presented: one aiming at coverage, a second at uniformity and a third at cardinality. All algorithms are based on the $\epsilon$-constraint method and are insensitive to the quality of nadir point estimation.

The algorithms were tested on 240 binary knapsack instances, 150 integer knapsack instances and 180 general problem instances, with three and four objectives, on the ability to efficiently obtain the full Pareto front. The algorithms that target uniformity and cardinality demonstrated to be very efficient and among the best in literature. The algorithms were also tested on the quality of the obtained representation, having demonstrated that, although the coverage and uniformity algorithms showed good results on their targeted measures, the cardinality algorithm is more flexible, privileging, depending on the cardinality, more the uniformity or the coverage of the representation.

As future work, we intend to study the incorporation of this representation strategies in interactive algorithms. Namely, we believe that the cardinality algorithm will provide a good basis in the sense that, in a first stage a uniform distribution over the whole Pareto front can be obtained and, at a later stage, as the DM's preferences are refined, a representation with a good coverage of the selected area can be obtained. At last, we would like to apply this representation algorithms in real multi-objective problems that would benefit from being treated as such, without the computational complexity of computing the whole Pareto front and difficulty of analysing it.

\section*{Acknowledgements}
\noindent Mariana Mesquita-Cunha acknowledges the support by national funds through FCT, under the research grant SFRH/BD/149441/2019. José Rui Figueira acknowledges the support by national funds through FCT, under DOME research project, PTDC/CCI-COM/31198/2017. Ana Paula Barbosa-Póvoa and Mariana Mesquita-Cunha acknowledge the support by national funds through FCT, under the Data2Help research project, DSAIPA/AI/0044/2018, and the project 1801P.00740 PTDC/EGE-OGE/28071/2017 - LISBOA-01-0145-FEDER-028071.
\bibliographystyle{elsarticle-harv}
\bibliography{references}

\clearpage



\appendix
\noindent 
\section{Figures of the illustrative example}\label{sec:Appendix_B}
\begin{center}
\begin{figure*}[!ht]
    \centering
    \begin{subfigure}[b]{0.45\textwidth}
        \scalebox{.9}{\begin{tikzpicture}
	\definecolor{darkblue}{RGB}{0,121,173}
	\pgfplotsset{
		/pgfplots/xlabel near ticks/.style={ /pgfplots/every axis x label/.style={ at={(ticklabel cs:0.5)},anchor=near ticklabel,},},
		/pgfplots/ylabel near ticks/.style={ /pgfplots/every axis y label/.style={ at={(ticklabel cs:0.5)},rotate=90,anchor=near ticklabel},
		}
	}
	\begin{axis}[%
	name=plot2,
	xmin=-38, xmax=59, ymin=-58, ymax=52, xtick={-28,0,49},
	ytick={-48,0,42},
	yticklabels={$-48$,$0$,$42$},
	xticklabels={$-28$, $0$, $49$}]

	\addplot +[red, no markers, name path=vazul] coordinates {(-38, -48) (59, -48)};
	\addplot [blue, no markers, name path=hverde] coordinates {(-28,-58) (-28,52)};
	\addplot[only marks,mark=*,mark size=2.9pt, color=darkblue, fill=darkblue] coordinates {(9,-14)};
	\end{axis}
\end{tikzpicture}}
        \caption{First iteration over objective $z_3$}    
        \label{fig:GPBA-A-1}
    \end{subfigure}%
    ~ 
    \begin{subfigure}[b]{0.45\textwidth}  
        \scalebox{.9}{\begin{tikzpicture}
	\definecolor{darkblue}{RGB}{0,121,173}
	\pgfplotsset{
		/pgfplots/xlabel near ticks/.style={ /pgfplots/every axis x label/.style={ at={(ticklabel cs:0.5)},anchor=near ticklabel,},},
		/pgfplots/ylabel near ticks/.style={ /pgfplots/every axis y label/.style={ at={(ticklabel cs:0.5)},rotate=90,anchor=near ticklabel},
			}
	}
	\begin{axis}[%
	name=plot2,
	xmin=-38, xmax=59, ymin=-58, ymax=52, xtick={-28,0,49},
	ytick={-48,0,42},
	yticklabels={$-48$,$0$,$42$},
	xticklabels={$-28$, $0$, $49$}]

	\addplot +[red, no markers, name path=vazul] coordinates {(-38, 42) (59, 42)};
	\addplot [blue, no markers, name path=hverde] coordinates {(-28,-58) (-28,52)};
	
	\addplot[only marks,mark=*,mark size=2.9pt, color=darkblue, fill=darkblue] coordinates {(9,-14)};
	\addplot[only marks,mark=*,mark size=2.9pt, color=darkblue, fill=darkblue] coordinates {(20,42)};
	\end{axis}
\end{tikzpicture}}
        \caption{Second iteration over objective $z_3$}
        \label{fig:GPBA-A-2}
    \end{subfigure}
    \begin{subfigure}[b]{0.45\textwidth}   
        \scalebox{.9}{\begin{tikzpicture}
	\definecolor{darkblue}{RGB}{0,121,173}
	\pgfplotsset{
		/pgfplots/xlabel near ticks/.style={ /pgfplots/every axis x label/.style={ at={(ticklabel cs:0.5)},anchor=near ticklabel,},},
		/pgfplots/ylabel near ticks/.style={ /pgfplots/every axis y label/.style={ at={(ticklabel cs:0.5)},rotate=90,anchor=near ticklabel},
			}
	}
	\begin{axis}[%
	name=plot2,
	xmin=-38, xmax=59, ymin=-58, ymax=52, xtick={-28,0,49},
	ytick={-48,0,42},
	yticklabels={$-48$,$0$,$42$},
	xticklabels={$-28$, $0$, $49$}]

	\addplot +[red, no markers, name path=vazul] coordinates {(-38, 14) (59, 14)};
	\addplot [blue, no markers, name path=hverde] coordinates {(-28,-58) (-28,52)};
	
	\addplot[only marks,mark=*,mark size=2.9pt, color=darkblue, fill=darkblue] coordinates {(9,-14)};
	\addplot[only marks,mark=*,mark size=2.9pt, color=darkblue, fill=darkblue] coordinates {(20,42)};
	\addplot[only marks,mark=*,mark size=2.9pt, color=darkblue, fill=darkblue] coordinates {(13,14)};
	\end{axis}
\end{tikzpicture}}
        \caption{Third iteration over objective $z_3$}    
        \label{fig:GPBA-A-3}
    \end{subfigure}%
    ~ 
    \begin{subfigure}[b]{0.45\textwidth}   
        \scalebox{.9}{\begin{tikzpicture}
	\definecolor{darkblue}{RGB}{0,121,173}
	\pgfplotsset{
		/pgfplots/xlabel near ticks/.style={ /pgfplots/every axis x label/.style={ at={(ticklabel cs:0.5)},anchor=near ticklabel,},},
		/pgfplots/ylabel near ticks/.style={ /pgfplots/every axis y label/.style={ at={(ticklabel cs:0.5)},rotate=90,anchor=near ticklabel},
			}
	}
	\begin{axis}[%
	name=plot2,
	xmin=-38, xmax=59, ymin=-58, ymax=52, xtick={-28,0,49},
	ytick={-48,0,42},
	yticklabels={$-48$,$0$,$42$},
	xticklabels={$-28$, $0$, $49$}]

	\addplot +[red, no markers, name path=vazul] coordinates {(-38, 0) (59, 0)};
	
	\addplot [blue, no markers, name path=hverde] coordinates {(-28,-58) (-28,52)};
	
	\addplot[only marks,mark=*,mark size=2.9pt, color=darkblue, fill=darkblue] coordinates {(9,-14)};
	\addplot[only marks,mark=*,mark size=2.9pt, color=darkblue, fill=darkblue] coordinates {(20,42)};
	\addplot[only marks,mark=*,mark size=2.9pt, color=darkblue, fill=darkblue] coordinates {(13,14)};
	\addplot[only marks,mark=*,mark size=2.9pt, color=darkblue, fill=darkblue] coordinates {(6,1)};
	\end{axis}
\end{tikzpicture}}
        \caption{Forth iteration over objective $z_3$}   
        \label{fig:GPBA-A-4}
    \end{subfigure}
    \begin{subfigure}[b]{0.45\textwidth}   
        \scalebox{.9}{\begin{tikzpicture}
	\definecolor{darkblue}{RGB}{0,121,173}
	\pgfplotsset{
		/pgfplots/xlabel near ticks/.style={ /pgfplots/every axis x label/.style={ at={(ticklabel cs:0.5)},anchor=near ticklabel,},},
		/pgfplots/ylabel near ticks/.style={ /pgfplots/every axis y label/.style={ at={(ticklabel cs:0.5)},rotate=90,anchor=near ticklabel},
			}
	}
	\begin{axis}[%
	name=plot2,
	xmin=-38, xmax=59, ymin=-58, ymax=52, xtick={-28,0,49},
	ytick={-48,0,42},
	yticklabels={$-48$,$0$,$42$},
	xticklabels={$-28$, $0$, $49$}]

	\addplot +[red, no markers, name path=vazul] coordinates {(-38, 28) (59, 28)};
	
	\addplot [blue, no markers, name path=hverde] coordinates {(-28,-58) (-28,52)};
	
	\addplot[only marks,mark=*,mark size=2.9pt, color=darkblue, fill=darkblue] coordinates {(9,-14)};
	\addplot[only marks,mark=*,mark size=2.9pt, color=darkblue, fill=darkblue] coordinates {(20,42)};
	\addplot[only marks,mark=*,mark size=2.9pt, color=darkblue, fill=darkblue] coordinates {(13,14)};
	\addplot[only marks,mark=*,mark size=2.9pt, color=darkblue, fill=darkblue] coordinates {(6,1)};
	\addplot[only marks,mark=*,mark size=2.9pt, color=darkblue, fill=darkblue] coordinates {(13,29)};
	\end{axis}
\end{tikzpicture}}
        \caption{Fifth iteration over objective $z_3$}   
        \label{fig:GPBA-A-5}
    \end{subfigure}
    \caption{Illustrative example for the coverage representation, \textit{GPBA-A}. The acceptable coverage error was defined as $\gamma=15$ in each objective. Only the first iteration over the outermost loop is depicted.}
    \label{fig:Example-GPBA-A}
\end{figure*}
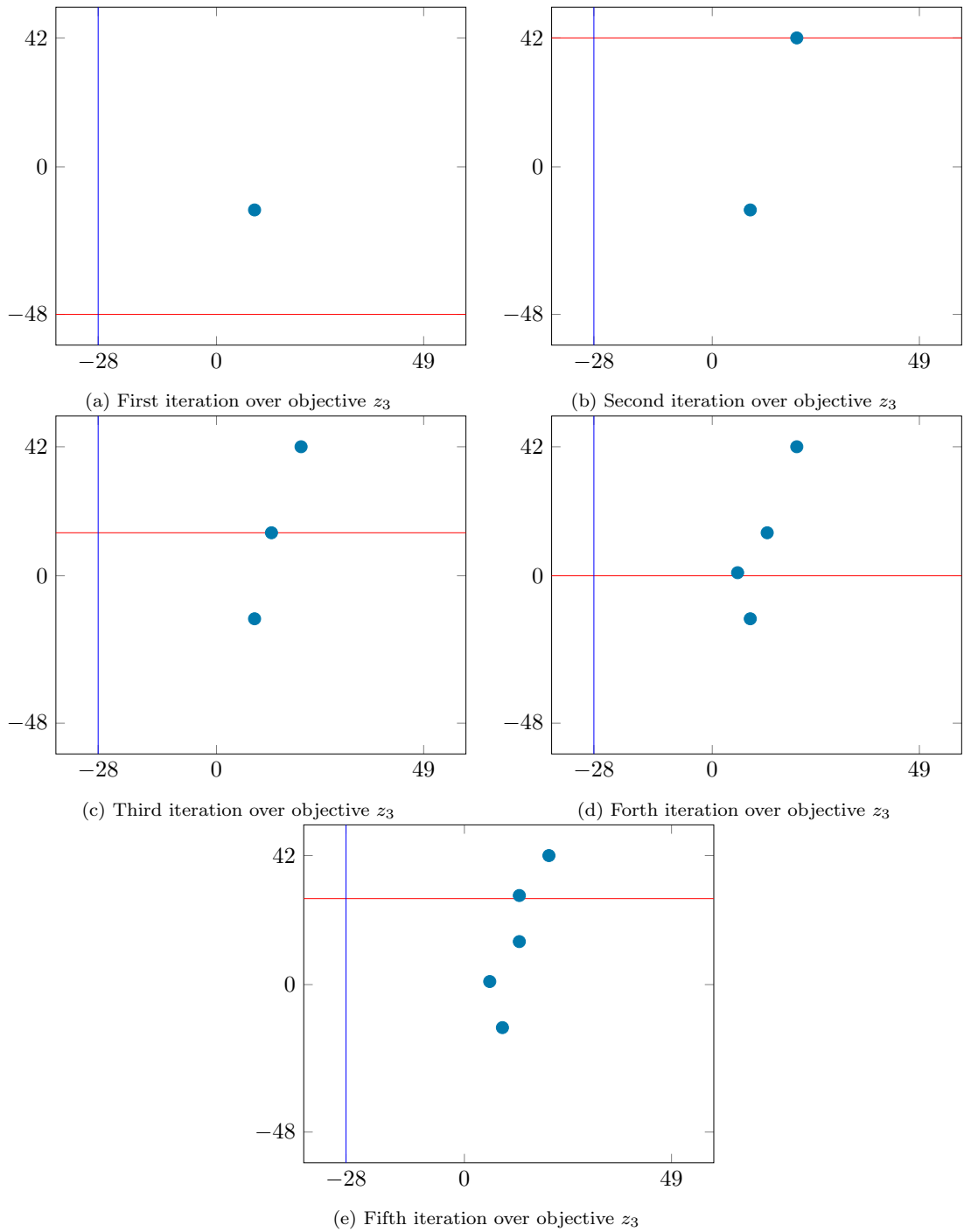
\end{center}

\begin{center}
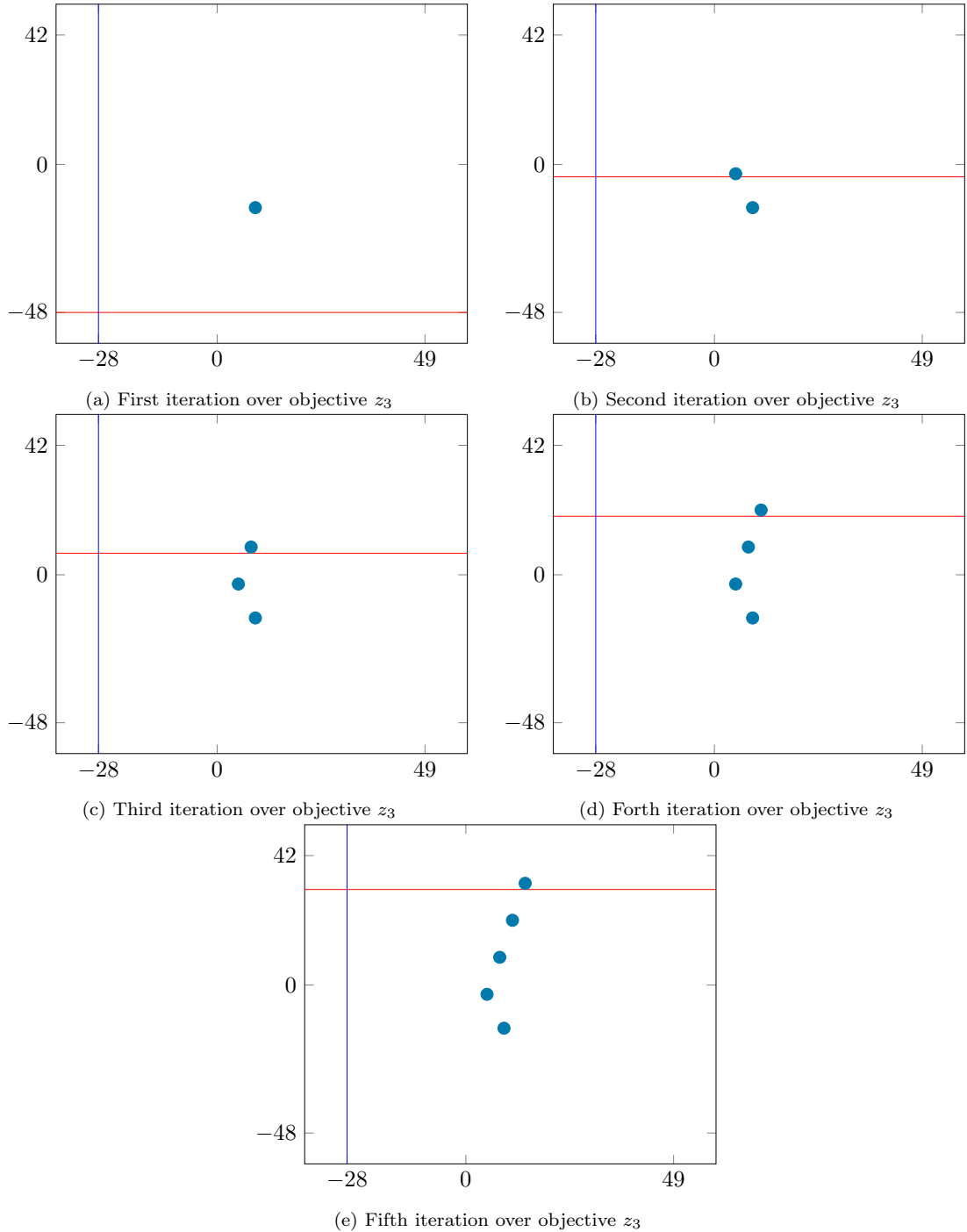
\begin{figure*}[!ht]
    \centering
    \begin{subfigure}[b]{0.45\textwidth}
        \scalebox{.9}{\begin{tikzpicture}
	\definecolor{darkblue}{RGB}{0,121,173}
	\pgfplotsset{
		/pgfplots/xlabel near ticks/.style={ /pgfplots/every axis x label/.style={ at={(ticklabel cs:0.5)},anchor=near ticklabel,},},
		/pgfplots/ylabel near ticks/.style={ /pgfplots/every axis y label/.style={ at={(ticklabel cs:0.5)},rotate=90,anchor=near ticklabel},
		}
	}
	\begin{axis}[%
	name=plot2,
	xmin=-38, xmax=59, ymin=-58, ymax=52, xtick={-28,0,49},
	ytick={-48,0,42},
	yticklabels={$-48$,$0$,$42$},
	xticklabels={$-28$, $0$, $49$}]

	\addplot +[red, no markers, name path=vazul] coordinates {(-38, -48) (59, -48)};
	\addplot [blue, no markers, name path=hverde] coordinates {(-28,-58) (-28,52)};
	\addplot[only marks,mark=*,mark size=2.9pt, color=darkblue, fill=darkblue] coordinates {(9,-14)};
	\end{axis}
\end{tikzpicture}}
        \caption{First iteration over objective $z_3$}    
        \label{fig:GPBA-B-1}
    \end{subfigure}%
    ~ 
    \begin{subfigure}[b]{0.45\textwidth}  
        \scalebox{.9}{\begin{tikzpicture}
	\definecolor{darkblue}{RGB}{0,121,173}
	\pgfplotsset{
		/pgfplots/xlabel near ticks/.style={ /pgfplots/every axis x label/.style={ at={(ticklabel cs:0.5)},anchor=near ticklabel,},},
		/pgfplots/ylabel near ticks/.style={ /pgfplots/every axis y label/.style={ at={(ticklabel cs:0.5)},rotate=90,anchor=near ticklabel},
		}
	}
	\begin{axis}[%
	name=plot2,
	xmin=-38, xmax=59, ymin=-58, ymax=52, xtick={-28,0,49},
	ytick={-48,0,42},
	yticklabels={$-48$,$0$,$42$},
	xticklabels={$-28$, $0$, $49$}]

	\addplot +[red, no markers, name path=vazul] coordinates {(-38, -4) (59, -4)};
	\addplot [blue, no markers, name path=hverde] coordinates {(-28,-58) (-28,52)};
	\addplot[only marks,mark=*,mark size=2.9pt, color=darkblue, fill=darkblue] coordinates {(9,-14)};
	\addplot[only marks,mark=*,mark size=2.9pt, color=darkblue, fill=darkblue] coordinates {(5,-3)};
	\end{axis}
\end{tikzpicture}}
        \caption{Second iteration over objective $z_3$}
        \label{fig:GPBA-B-2}
    \end{subfigure}
    \begin{subfigure}[b]{0.45\textwidth}   
        \scalebox{.9}{\begin{tikzpicture}
	\definecolor{darkblue}{RGB}{0,121,173}
	\pgfplotsset{
		/pgfplots/xlabel near ticks/.style={ /pgfplots/every axis x label/.style={ at={(ticklabel cs:0.5)},anchor=near ticklabel,},},
		/pgfplots/ylabel near ticks/.style={ /pgfplots/every axis y label/.style={ at={(ticklabel cs:0.5)},rotate=90,anchor=near ticklabel},
		}
	}
	\begin{axis}[%
	name=plot2,
	xmin=-38, xmax=59, ymin=-58, ymax=52, xtick={-28,0,49},
	ytick={-48,0,42},
	yticklabels={$-48$,$0$,$42$},
	xticklabels={$-28$, $0$, $49$}]

	\addplot +[red, no markers, name path=vazul] coordinates {(-38, 7) (59, 7)};
	\addplot [blue, no markers, name path=hverde] coordinates {(-28,-58) (-28,52)};
	\addplot[only marks,mark=*,mark size=2.9pt, color=darkblue, fill=darkblue] coordinates {(9,-14)};
	\addplot[only marks,mark=*,mark size=2.9pt, color=darkblue, fill=darkblue] coordinates {(5,-3)};
	\addplot[only marks,mark=*,mark size=2.9pt, color=darkblue, fill=darkblue] coordinates {(8,9)};
	\end{axis}
\end{tikzpicture}}
        \caption{Third iteration over objective $z_3$}    
        \label{fig:GPBA-B-3}
    \end{subfigure}%
    ~ 
    \begin{subfigure}[b]{0.45\textwidth}   
        \scalebox{.9}{\begin{tikzpicture}
	\definecolor{darkblue}{RGB}{0,121,173}
	\pgfplotsset{
		/pgfplots/xlabel near ticks/.style={ /pgfplots/every axis x label/.style={ at={(ticklabel cs:0.5)},anchor=near ticklabel,},},
		/pgfplots/ylabel near ticks/.style={ /pgfplots/every axis y label/.style={ at={(ticklabel cs:0.5)},rotate=90,anchor=near ticklabel},
		}
	}
	\begin{axis}[%
	name=plot2,
	xmin=-38, xmax=59, ymin=-58, ymax=52, xtick={-28,0,49},
	ytick={-48,0,42},
	yticklabels={$-48$,$0$,$42$},
	xticklabels={$-28$, $0$, $49$}]

	\addplot +[red, no markers, name path=vazul] coordinates {(-38, 19) (59, 19)};
	\addplot [blue, no markers, name path=hverde] coordinates {(-28,-58) (-28,52)};
	\addplot[only marks,mark=*,mark size=2.9pt, color=darkblue, fill=darkblue] coordinates {(9,-14)};
	\addplot[only marks,mark=*,mark size=2.9pt, color=darkblue, fill=darkblue] coordinates {(5,-3)};
	\addplot[only marks,mark=*,mark size=2.9pt, color=darkblue, fill=darkblue] coordinates {(8,9)};
	\addplot[only marks,mark=*,mark size=2.9pt, color=darkblue, fill=darkblue] coordinates {(11,21)};
	\end{axis}
\end{tikzpicture}}
        \caption{Forth iteration over objective $z_3$}   
        \label{fig:GPBA-B-4}
    \end{subfigure}
    \begin{subfigure}[b]{0.45\textwidth}   
        \scalebox{.9}{\begin{tikzpicture}
	\definecolor{darkblue}{RGB}{0,121,173}
	\pgfplotsset{
		/pgfplots/xlabel near ticks/.style={ /pgfplots/every axis x label/.style={ at={(ticklabel cs:0.5)},anchor=near ticklabel,},},
		/pgfplots/ylabel near ticks/.style={ /pgfplots/every axis y label/.style={ at={(ticklabel cs:0.5)},rotate=90,anchor=near ticklabel},
		}
	}
	\begin{axis}[%
	name=plot2,
	xmin=-38, xmax=59, ymin=-58, ymax=52, xtick={-28,0,49},
	ytick={-48,0,42},
	yticklabels={$-48$,$0$,$42$},
	xticklabels={$-28$, $0$, $49$}]

	\addplot +[red, no markers, name path=vazul] coordinates {(-38, 31) (59, 31)};
	\addplot [blue, no markers, name path=hverde] coordinates {(-28,-58) (-28,52)};
	\addplot[only marks,mark=*,mark size=2.9pt, color=darkblue, fill=darkblue] coordinates {(9,-14)};
	\addplot[only marks,mark=*,mark size=2.9pt, color=darkblue, fill=darkblue] coordinates {(5,-3)};
	\addplot[only marks,mark=*,mark size=2.9pt, color=darkblue, fill=darkblue] coordinates {(8,9)};
	\addplot[only marks,mark=*,mark size=2.9pt, color=darkblue, fill=darkblue] coordinates {(11,21)};
	\addplot[only marks,mark=*,mark size=2.9pt, color=darkblue, fill=darkblue] coordinates {(14,33)};
	\end{axis}
\end{tikzpicture}}
        \caption{Fifth iteration over objective $z_3$}   
        \label{fig:GPBA-B-5}
    \end{subfigure}
    \caption{Illustrative example for the uniformity representation, \textit{GPBA-B}. The uniformity level defined was $\delta=10$ in all objectives. Only the first iteration over the outermost loop is depicted.}
    \label{fig:Example-GPBA-B}
\end{figure*}
\end{center}

\begin{center}
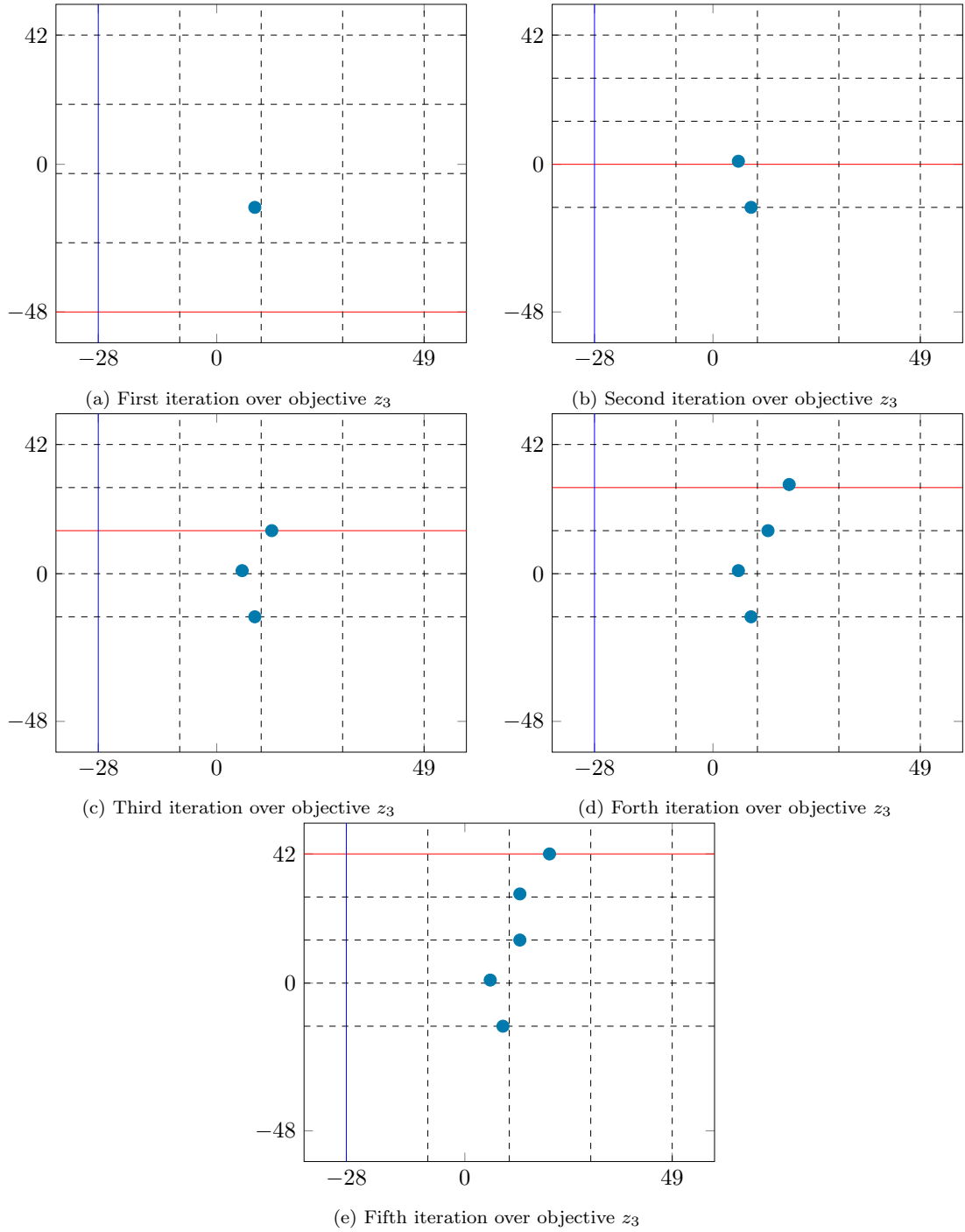
\begin{figure*}[!ht]
    \centering
    \begin{subfigure}[b]{0.45\textwidth}
        \scalebox{.9}{\begin{tikzpicture}
	\definecolor{darkblue}{RGB}{0,121,173}
	\pgfplotsset{
		/pgfplots/xlabel near ticks/.style={ /pgfplots/every axis x label/.style={ at={(ticklabel cs:0.5)},anchor=near ticklabel,},},
		/pgfplots/ylabel near ticks/.style={ /pgfplots/every axis y label/.style={ at={(ticklabel cs:0.5)},rotate=90,anchor=near ticklabel},
		}
	}
	\begin{axis}[%
	name=plot2,
	xmin=-38, xmax=59, ymin=-58, ymax=52, xtick={-28,0,49},
	ytick={-48,0,42},
	yticklabels={$-48$,$0$,$42$},
	xticklabels={$-28$, $0$, $49$}]

	\addplot +[red, no markers, name path=vazul] coordinates {(-38, -48) (59, -48)};
	\addplot [blue, no markers, name path=hverde] coordinates {(-28,-58) (-28,52)};
	
	\addplot [dashed, no markers] coordinates {(-8.75,-58) (-8.75,52)};
	\addplot [dashed, no markers] coordinates {(10.5,-58) (10.5,52)};
	\addplot [dashed, no markers] coordinates {(29.75,-58) (29.75,52)};
	\addplot [dashed, no markers] coordinates {(49,-58) (49,52)};
	\addplot [dashed, no markers] coordinates {(-38, -25.5) (59, -25.5)};
	\addplot [dashed, no markers] coordinates {(-38, -3) (59, -3)};
	\addplot [dashed, no markers] coordinates {(-38, 19.5) (59, 19.5)};
	\addplot [dashed, no markers] coordinates {(-38, 42) (59, 42)};
	
	\addplot[only marks,mark=*,mark size=2.9pt, color=darkblue, fill=darkblue] coordinates {(9,-14)};
	\end{axis}
\end{tikzpicture}}
        \caption{First iteration over objective $z_3$}    
        \label{fig:GPBA-C-1}
    \end{subfigure}%
    ~ 
    \begin{subfigure}[b]{0.45\textwidth}  
        \scalebox{.9}{\begin{tikzpicture}
	\definecolor{darkblue}{RGB}{0,121,173}
	\pgfplotsset{
		/pgfplots/xlabel near ticks/.style={ /pgfplots/every axis x label/.style={ at={(ticklabel cs:0.5)},anchor=near ticklabel,},},
		/pgfplots/ylabel near ticks/.style={ /pgfplots/every axis y label/.style={ at={(ticklabel cs:0.5)},rotate=90,anchor=near ticklabel},
		}
	}
	\begin{axis}[%
	name=plot2,
	xmin=-38, xmax=59, ymin=-58, ymax=52, xtick={-28,0,49},
	ytick={-48,0,42},
	yticklabels={$-48$,$0$,$42$},
	xticklabels={$-28$, $0$, $49$}]

	\addplot +[red, no markers, name path=vazul] coordinates {(-38, 0) (59, 0)};
	\addplot [blue, no markers, name path=hverde] coordinates {(-28,-58) (-28,52)};
	
	\addplot [dashed, no markers] coordinates {(-8.75,-58) (-8.75,52)};
	\addplot [dashed, no markers] coordinates {(10.5,-58) (10.5,52)};
	\addplot [dashed, no markers] coordinates {(29.75,-58) (29.75,52)};
	\addplot [dashed, no markers] coordinates {(49,-58) (49,52)};
	\addplot [dashed, no markers] coordinates {(-38, -14) (59, -14)};
	\addplot [dashed, no markers] coordinates {(-38, 14) (59, 14)};
	\addplot [dashed, no markers] coordinates {(-38, 28) (59, 28)};
	\addplot [dashed, no markers] coordinates {(-38, 42) (59, 42)};
	
	\addplot[only marks,mark=*,mark size=2.9pt, color=darkblue, fill=darkblue] coordinates {(9,-14)};
	\addplot[only marks,mark=*,mark size=2.9pt, color=darkblue, fill=darkblue] coordinates {(6,1)};
	\end{axis}
\end{tikzpicture}}
        \caption{Second iteration over objective $z_3$}
        \label{fig:GPBA-C-2}
    \end{subfigure}
    \begin{subfigure}[b]{0.45\textwidth}   
        \scalebox{.9}{\begin{tikzpicture}
	\definecolor{darkblue}{RGB}{0,121,173}
	\pgfplotsset{
		/pgfplots/xlabel near ticks/.style={ /pgfplots/every axis x label/.style={ at={(ticklabel cs:0.5)},anchor=near ticklabel,},},
		/pgfplots/ylabel near ticks/.style={ /pgfplots/every axis y label/.style={ at={(ticklabel cs:0.5)},rotate=90,anchor=near ticklabel},
		}
	}
	\begin{axis}[%
	name=plot2,
	xmin=-38, xmax=59, ymin=-58, ymax=52, xtick={-28,0,49},
	ytick={-48,0,42},
	yticklabels={$-48$,$0$,$42$},
	xticklabels={$-28$, $0$, $49$}]

	\addplot +[red, no markers, name path=vazul] coordinates {(-38, 14) (59, 14)};
	\addplot [blue, no markers, name path=hverde] coordinates {(-28,-58) (-28,52)};
	
	\addplot [dashed, no markers] coordinates {(-8.75,-58) (-8.75,52)};
	\addplot [dashed, no markers] coordinates {(10.5,-58) (10.5,52)};
	\addplot [dashed, no markers] coordinates {(29.75,-58) (29.75,52)};
	\addplot [dashed, no markers] coordinates {(49,-58) (49,52)};
	\addplot [dashed, no markers] coordinates {(-38, -14) (59, -14)};
	\addplot [dashed, no markers] coordinates {(-38, 0) (59, 0)};
	\addplot [dashed, no markers] coordinates {(-38, 28) (59, 28)};
	\addplot [dashed, no markers] coordinates {(-38, 42) (59, 42)};
	
	\addplot[only marks,mark=*,mark size=2.9pt, color=darkblue, fill=darkblue] coordinates {(9,-14)};
	\addplot[only marks,mark=*,mark size=2.9pt, color=darkblue, fill=darkblue] coordinates {(6,1)};
	\addplot[only marks,mark=*,mark size=2.9pt, color=darkblue, fill=darkblue] coordinates {(13,14)};
	\end{axis}
\end{tikzpicture}}
        \caption{Third iteration over objective $z_3$}    
        \label{fig:GPBA-C-3}
    \end{subfigure}%
    ~ 
    \begin{subfigure}[b]{0.45\textwidth}   
        \scalebox{.9}{\begin{tikzpicture}
	\definecolor{darkblue}{RGB}{0,121,173}
	\pgfplotsset{
		/pgfplots/xlabel near ticks/.style={ /pgfplots/every axis x label/.style={ at={(ticklabel cs:0.5)},anchor=near ticklabel,},},
		/pgfplots/ylabel near ticks/.style={ /pgfplots/every axis y label/.style={ at={(ticklabel cs:0.5)},rotate=90,anchor=near ticklabel},
		}
	}
	\begin{axis}[%
	name=plot2,
	xmin=-38, xmax=59, ymin=-58, ymax=52, xtick={-28,0,49},
	ytick={-48,0,42},
	yticklabels={$-48$,$0$,$42$},
	xticklabels={$-28$, $0$, $49$}]

	\addplot +[red, no markers, name path=vazul] coordinates {(-38, 28) (59, 28)};
	\addplot [blue, no markers, name path=hverde] coordinates {(-28,-58) (-28,52)};
	
	\addplot [dashed, no markers] coordinates {(-8.75,-58) (-8.75,52)};
	\addplot [dashed, no markers] coordinates {(10.5,-58) (10.5,52)};
	\addplot [dashed, no markers] coordinates {(29.75,-58) (29.75,52)};
	\addplot [dashed, no markers] coordinates {(49,-58) (49,52)};
	\addplot [dashed, no markers] coordinates {(-38, -14) (59, -14)};
	\addplot [dashed, no markers] coordinates {(-38, 0) (59, 0)};
	\addplot [dashed, no markers] coordinates {(-38, 14) (59, 14)};
	\addplot [dashed, no markers] coordinates {(-38, 42) (59, 42)};
	
	\addplot[only marks,mark=*,mark size=2.9pt, color=darkblue, fill=darkblue] coordinates {(9,-14)};
	\addplot[only marks,mark=*,mark size=2.9pt, color=darkblue, fill=darkblue] coordinates {(6,1)};
	\addplot[only marks,mark=*,mark size=2.9pt, color=darkblue, fill=darkblue] coordinates {(13,14)};
	\addplot[only marks,mark=*,mark size=2.9pt, color=darkblue, fill=darkblue] coordinates {(18,29)};
	\end{axis}
\end{tikzpicture}}
        \caption{Forth iteration over objective $z_3$}   
        \label{fig:GPBA-C-4}
    \end{subfigure}
    \begin{subfigure}[b]{0.45\textwidth}   
        \scalebox{.9}{\begin{tikzpicture}
	\definecolor{darkblue}{RGB}{0,121,173}
	\pgfplotsset{
		/pgfplots/xlabel near ticks/.style={ /pgfplots/every axis x label/.style={ at={(ticklabel cs:0.5)},anchor=near ticklabel,},},
		/pgfplots/ylabel near ticks/.style={ /pgfplots/every axis y label/.style={ at={(ticklabel cs:0.5)},rotate=90,anchor=near ticklabel},
		}
	}
	\begin{axis}[%
	name=plot2,
	xmin=-38, xmax=59, ymin=-58, ymax=52, xtick={-28,0,49},
	ytick={-48,0,42},
	yticklabels={$-48$,$0$,$42$},
	xticklabels={$-28$, $0$, $49$}]

	\addplot +[red, no markers, name path=vazul] coordinates {(-38, 42) (59, 42)};
	\addplot [blue, no markers, name path=hverde] coordinates {(-28,-58) (-28,52)};
	
	\addplot [dashed, no markers] coordinates {(-8.75,-58) (-8.75,52)};
	\addplot [dashed, no markers] coordinates {(10.5,-58) (10.5,52)};
	\addplot [dashed, no markers] coordinates {(29.75,-58) (29.75,52)};
	\addplot [dashed, no markers] coordinates {(49,-58) (49,52)};
	\addplot [dashed, no markers] coordinates {(-38, -14) (59, -14)};
	\addplot [dashed, no markers] coordinates {(-38, 0) (59, 0)};
	\addplot [dashed, no markers] coordinates {(-38, 14) (59, 14)};
	\addplot [dashed, no markers] coordinates {(-38, 28) (59, 28)};
	
	\addplot[only marks,mark=*,mark size=2.9pt, color=darkblue, fill=darkblue] coordinates {(9,-14)};
	\addplot[only marks,mark=*,mark size=2.9pt, color=darkblue, fill=darkblue] coordinates {(6,1)};
	\addplot[only marks,mark=*,mark size=2.9pt, color=darkblue, fill=darkblue] coordinates {(13,14)};
	\addplot[only marks,mark=*,mark size=2.9pt, color=darkblue, fill=darkblue] coordinates {(13,29)};
	\addplot[only marks,mark=*,mark size=2.9pt, color=darkblue, fill=darkblue] coordinates {(20,42)};
	\end{axis}
\end{tikzpicture}}
        \caption{Fifth iteration over objective $z_3$}   
        \label{fig:GPBA-C-5}
    \end{subfigure}
    \caption{Illustrative example for the cardinality representation, \textit{GPBA-C}. The desired cardinality defined was $5$ in all objectives. Only the first iteration over the outermost loop is depicted.}
    \label{fig:Example-GPBA-C}
\end{figure*}
\end{center}

\end{document}